\magnification=1200
\overfullrule 0cm

\def\ref#1{\lbrack {#1}\rbrack}

\def\ekv#1#2{$${#2}\eqno(#1)$$}
\def\eekv#1#2#3{$$\eqalignno{&{#2}&({#1})\cr &{#3}\cr}$$}
\def\eeekv#1#2#3#4{$$\eqalignno{&{#2}&({#1})\cr &{#3}\cr &{#4}\cr}$$}

\def\iint{\int\hskip -2mm\int}

\font\liten=cmr10 at 8pt
\font\stor=cmr10 at 12pt

\def\aby{arbitrary}
\def\ably{arbitrarily}
\def\asy{asymptotic}
\def\bdd{bounded}
\def\bdy{boundary}

\def\canform{canonical transformation}

\def\diff{diffeomorphism}
\def\diffeo{diffeomorphism}

\def\ev{eigen-value}
\def\e{equation}
\def\fy{family}
\def\fu{function}

\def\fop{Fourier integral operator}

\def\hol{holomorphic}
\def\indep{independent}
\def\lhs{left hand side}
\def\mfld{manifold}
\def\neigh{neighborhood}
\def\nondeg{non-degenerate}
\def\op{operator}
\def\og{orthogonal}
\def\pb{problem}

\def\prop{proposition}
\def\Prop{Proposition}
\def\pol{polynomial}
\def\pop{pseudodifferential operator}

\def\res{resonance}

\def\sa{self-adjoint}

\def\strans{^\sigma \hskip -2pt}

\def\sufly{sufficiently}
\def\tf{transformation}
\def\Th{Theorem}
\def\th{theorem}
\def\tf{transform}
\def\trans{^t\hskip -2pt}

\def\vf{vector field}
\def\wrt{with respect to}

\def\Re{{\rm Re\,}}
\def\Im{{\rm Im\,}}

\def\In{0}
\def\Sy{1}
\def\Eq{2}
\def\LF{3}
\def\BN{4}
\def\Pa{5}

\centerline{\stor Birkhoff normal forms for \fop{}s II.}
\medskip
\smallskip
\centerline{\bf A. Iantchenko\footnote{*}{\rm Malm{\"o} H{\"o}gskola, Teknik
och Samh{\"a}lle, SE-20506 Malm{\"o}}
\rm and \bf J. Sj{\"o}strand\footnote{$\dagger$}{\rm Centre de
Math{\'e}matiques, Ecole Polytechnique, FR-91128 Palaiseau Cedex,\break and
UMR 7640 of CNRS}\footnote{}{\it Key words: \liten \fop{}, Birkhoff,
normal form, symplectic.}\footnote{}{\it MSC2000: \liten 34C20, 34K17,
35S05, 35S30, 37J40, 70K45, 81Q20}}
\bigskip
\par\noindent \it Abstract. \liten In this work we construct
logarithms and Birkhoff normal forms for elliptic \fop{}s in the
semi-classical limit, under more general assumptions than in a
previous work by the first author. The methods are similar but
slightly different.
\medskip
\par\noindent \it R{\'e}sum{\'e}. \liten Dans ce travail nous construisons
des logarithmes et des formes normales de Birkhoff pour des op{\'e}rateurs
int{\'e}graux de Fourier dans la limite semi-classique. Les hypoth{\`e}ses
sont plus faibles que dans un travail ant{\'e}rieur du premier auteur
et les m{\'e}thodes sont semblables mais un peu diff{\'e}rentes.\rm
\bigskip
\par
\centerline{\bf \In{}. Introduction.} 
\medskip
\par This work is a continuation of the work [Ia] of the first author.
As in that paper we consider the problem of constructing a microlocal
logarithm and a normal form of an elliptic semi-classical \fop{} near a
fixed point of the corresponding \canform{}. In [Ia] the \canform{} was
assumed to be of real hyperbolic type and the purpose of the present
work is to relax this assumption considerably, to what we think are the
natural conditions.

\par As in [Ia] a motivation is to improve the analysis of quantized
billiard ball maps near closed trajectories for \bdy{} value problems.
Then the associated \canform{} is the corresponding classical
Poincar{\'e} map. In the introduction of [Ia] one such problem in the case
of scattering by obstacles was mentioned, where this \canform{} is of
hyperbolic type near its fixed point, but in many other cases
the \canform can be more \aby{}. (Such an improvement of spectral
asymptotics for non-degenerate potential wells was obtained in [Sj] by
means of a Birkhoff normal form.)

\par Another motivation for returning to this study comes from recent
works on inverse spectral problems by Guillemin [Gui] and Zelditch
[Ze1--3]. In these works, quantum normal forms of the whole (Laplace)
operator were constructed in a \neigh{} of the closed trajectory. We
believe that it may often be sufficient and technically easier to work
with normal forms of a corresponding monodromy operator, which coincides
with the quantized billiard ball map in the case of \bdy{} problems. For
general semi-classical problems, this operator was recently introduced
and studied in a more systematic fashion in [SjZw]. In particular one
might arrive at very nice trace formulae by combining the basic trace
formula of that paper with the normal form obtained in the present work.

\par It turned out to be somewhat difficult to extend the whole method of
[Ia] to the more general situation here. In [Ia] a major idea was to
use that the symplectic group is connected and work with deformations
from the the identity transformation to the given \canform{}. In the
hyperbolic case this can be done in a such a way that the intermediate
\tf{}s satisfy the assumptions for having a normal form. In the
general case considered here, we encounter singular values for the
deformation parameter where the conditions are not fulfilled. It seems
possible to circumvent this difficulty by complexifying the deformation
parameter and use corresponding slightly complex \tf{}s. Eventually
however we found a method allowing us to avoid deforming the differential
of the \canform{} but only the higher order part in its Taylor expansion.
In this way we have families of objects which satisfy our assumptions
everywhere.

\par Consider a semiclassical \fop{} $A$ of order 0, with an associated
\canform{}
$\kappa :{\rm neigh\,}((0,0),{\bf R}^{2n})\to {\rm neigh\,}((0,0),{\bf
R}^{2n})$ having (0,0) as a fixed point. Assume $A$ is elliptic at
(0,0). The set of \ev{}s of $d\kappa (0)$ is then closed under inversion
$\lambda \to 1/\lambda $ and under complex conjugation. Assume that to
the distinct $\lambda $ in the spectrum of $d\kappa (0,0)$, we can
associate a logartithm $\mu =\log
\lambda $ in such a way that inversion and complex conjugation
correspond to the map $\mu \to -\mu $ and to complex conjugation
respectiveley. (Notice that this assumption excludes the existence of
negative \ev{}s of $d\kappa (0,0)$.) Assume also that
\ekv{\In{}.1}
{
\sum k_j\mu _j=2\pi in,\ k_j,n\in{\bf Z}\Rightarrow \sum k_j\mu _j=0
}
Then a real version the Lewis--Sternberg theorem (see [St], [Fr] and
\Th{} \Sy{}.3 below) tells us that we can write $\kappa (\rho )=\exp
H_p(\rho )+{\cal O}(\rho ^\infty )$ ($\rho =(x,\xi )\in{\bf R}^{2n})$ for
some smooth and real-valued
$p={\cal O}(\rho ^2)$. The first result of this work says that under the
same assumptions, we can write $A\equiv e^{-iP/h}$ modulo an operator
which vanishes to infinite order at $(x,\xi )=(0,0)$, $h=0$, where $h>0$
is the small semi-classical parameter (\Th{} \LF{}.2). Here $P$ is a
semi-classical \pop{} of order 0 with $p$
as its leading symbol.

\par The second result is a  straight forward extension of the normal
form in [Sj] and says that under the non-\res{} condition (\BN{}.2)
below, the \op{} $P$ has a simple normal form in terms of certain action
\op{}s. (In (\In{}.1) we do not have to be very precise concerning the
enumeration of the logaritms of the \ev{}s of the linearization as long
as we have have one of each, modulo the sign, while in (\BN{}.2) we
enumerate one rather specific half of these logs. Also notice that the
combination of the two conditions takes the simple form (\BN{}.3).) We do
not expect that the exclusion of negative \ev{}s is a serious
restriction, for if such
\ev{}s are present, we can find specially adapted and explicit metaplectic
\op{}s $M$ and apply our results to $MA$.

\par Both our results were obtained by the first author ([Ia]) under the
assumption that $d\kappa (0,0)$ is of real-hyperbolic type.

\par For the reader's convenience we have taken the pain review some
well-known linear symplectic geometry in section \Sy{}, and in that
section we also give a (probably well-known) proof of the real
version of the Lewis-Sternberg theorem, which has a structure that
we can follow in the proof of the corresponding quantum result for
logarithms of \fop{}s.

\par It is beyond the scope of this work to consider convergence
questions related to the perturbation series that appear in connection
with normal forms. See for instance H. Ito [It] and G. Popov [Po1,2].

\par The plan of the paper is the following:
\smallskip
\par\noindent In section \Sy{} we review some standard facts about
linear symplectic geometry and add a few remarks for treating the real
case. We also review a proof of a real version of the Lewis--Sternberg
theorem, that we can later use as  guideline for the proof of the
quantum result.
\smallskip
\par\noindent In section \Eq{}, we introduce some notions of
equivalence that are used in the formulation of the main results. These
notions are essentially the same as in [Ia].
\smallskip
\par\noindent In section \LF{} we establish the main result about
logarithms of \fop{}s.
\smallskip
\par\noindent In section \BN{} we give the "Birkhoff" normal form for
the logarithm i.e. for a certain $h$-\pop{} of order 0.
\smallskip
\par\noindent In section \Pa{} we extend the results to the parameter
dependent case. In many genuinely semi-classical problems we do not
have any homogeneity, inferring that the Poincar{\'e} map is essentially
energy \indep{}, and then we cannot expect in general that our
arithmetic condions be fulfilled for all energies in some interval.
Consequently it is of interest to know that the results are valid to
infinite order at points where the conditions are fulfilled.

\bigskip

\centerline{\bf \Sy{}. Review of some symplectic
geometry.} 
\medskip

\par In this section, we review some more or less well-known
arguments that will later be extended to the quantized case. See [MeHa].

\par We start with the linear case and let $A:{\bf C}^{2n}\to {\bf
C}^{2n}$ be linear and symplectic in the sense that
\ekv{\Sy{}.1}
{
\sigma (Ax,Ay)=\sigma (x,y),\ x,y\in{\bf C}^{2n},
}
where $\sigma $ is the standard symplectic 2-form on ${\bf C}^{2n}$. We
recall that this implies that $\det A=1$, since $A$ will conserve the
volume form $\sigma ^n/n!$. When $n=1$ the converse is also true.

\par Let $E_\lambda ={\cal N}((A-\lambda )^{2n})={\rm Ker\,}((A-\lambda
)^{2n})$ be the generalized eigen-space associated to $\lambda \in{\bf
C}\setminus\{ 0\}$.
\medskip
\par\noindent
\bf\Prop{} \Sy{}.1. \it If $\lambda \mu \ne 1$, $\lambda ,\mu \in{\bf
C}\setminus\{ 0\}$, then $E_\lambda $ and $E_\mu $ are symplectically
orthogonal: $E_\lambda \perp^\sigma  E_\mu $.
\rm\medskip
\par\noindent \bf Proof. \rm Let $E_\lambda ^{(j)}=E_\lambda \cap {\cal
N}((A-\lambda )^j)$, so that
$$0\ne E_\lambda ^{(1)}\subset E_\lambda ^{(2)}\subset ..\subset
E_\lambda ^{(2n)}=E_\lambda .$$
Define $E_\mu ^{(k)}$ in the same way. Then for $x\in E_\lambda
^{(1)}$, $y\in E_\mu ^{(1)}$:
$$\sigma (x,y)=\sigma ({1\over \lambda }Ax,{1\over \mu }Ay)={1\over
\lambda \mu }\sigma (Ax,Ay)={1\over \lambda \mu }\sigma (x,y),$$
and since $1/\lambda \mu \ne 1$, we get $\sigma (x,y)=0$.

\par Assume that we have for some $m\ge 2$:
\ekv{{\rm P}_m}
{
\sigma (x,y)=0,\hbox{ for }x\in E_\lambda ^{(j)},\, y\in E_\mu ^{(k)},\
j+k\le m. }
We have just established (${\rm P}_2$).

\par Let $x\in E_\lambda ^{(j)}$, $y\in E_\mu ^{(k)}$, $j+k=m+1$. Write
${A_\vert}_{E_\lambda }=\lambda +N$,
${A_\vert}_{E_\mu }=\mu +M$, where $N,M$ are nilpotent with
$N:E_\lambda ^{(j)}\to E_\lambda ^{(j-1)}$, $M:E_\mu ^{(k)}\to E_\mu
^{(k-1)}$. Then
$$\eqalign{\sigma (x,y)&=\sigma ({1\over \lambda }Ax-{1\over \lambda
}Nx,{1\over \mu }Ay-{1\over \mu }My)\cr
&={1\over \lambda \mu }\sigma (x,y)+\sigma (-{N\over \lambda }x,{1\over
\mu} Ay)+\sigma ({1\over \lambda }Ax,-{1\over \mu }My)+\sigma ({1\over
\lambda }Nx,{1\over \mu }My).}$$
The last three terms vanish because of the induction assumption, and we
get $\sigma (x,y)=0$, so we have proved (${\rm P}_{m+1})$. The
proposition follows.\hfill{$\#$}\medskip

\par We conclude that $E_\lambda $ are isotropic if $\lambda ^2\ne 1$
(i.e. $\sigma $ vanishes on $E_\lambda \times E_\lambda $). We also see
that
$$E_1\perp^\sigma E_\lambda \hbox{ for }\lambda \ne 1,\
E_{-1}\perp^\sigma E_\mu \hbox{ for }\mu \ne -1. $$
It follows that
\ekv{\Sy{}.2}
{
{\bf C}^{2n}=E_1\oplus E_{-1}\oplus \bigoplus_1^k (E_{\lambda _j}\oplus
E_{1\over \lambda _j}), }
where $\lambda _j, \lambda _j^{-1}$ and possibly $1,-1$ denote the
distinct \ev{}s of $A$ with $\lambda _j\ne \pm 1$. Moreover, all the
$\oplus$ and
$\bigoplus$ except the $\oplus$s inside the parenthesies indicate
symplectically
\og{} decomposition.

\par For $\lambda \not\in \{ 1,-1\}$, write ${A_\vert}_{E_\lambda
}=\lambda +N_\lambda $ with $N_\lambda $ nilpotent. ${\sigma
_\vert}_{E_\lambda \times E_{1\over \lambda }}$ is \nondeg{} and if
${\strans A}$ denotes the symplectic transpose of $A$, so that $\sigma
(Ax,y)=\sigma (x,{\strans A}y)$, then ${\strans}{\strans A}=A$, and
${\strans A} A=1$, and hence
$${\strans A}=A^{-1}.$$
Also notice that the $E_\lambda $ are invariant subspaces for
$\strans{A}$. On $E_{1/\lambda }$ we have on the one hand
$A=1/\lambda +N_{1/\lambda }$ and on the other hand,
$$1={\strans A}A=(\lambda +{\strans N}_\lambda )({1\over \lambda
}+N_{{1\over \lambda }})=(1+{1\over \lambda }{\strans N}_\lambda
)(1+\lambda N_{1\over \lambda
}).$$
Hence,
\ekv{\Sy{}.3}
{\lambda N_{{1\over \lambda }}=-({1\over \lambda }{\strans N_\lambda
})+({1\over \lambda }{\strans N}_\lambda )^2-(..)^3+...\,\, .}

\par On $E_1$ we have $A=1+N_1$, where $N_1$ is nilpotent. The
requirement that $A$ be symplectic means that
\ekv{\Sy{}.4}
{(1+{\strans N_1})(1+N_1)=1.}
Similarly on $E_{-1}$, we have
\ekv{\Sy{}.5}
{({\strans N_{-1}}-1)(N_{-1}-1)=1.}

\par Conversely consider a decomposition of ${\bf C}^{2n}$ as in (\Sy{}.2)
into symplectically \og{} spaces $E_1$, $E_{-1}$, $E_{\lambda _j}\oplus
E_{1/\lambda _j}$ with $E_{\pm 1}$ symplectic, $E_{\lambda _j}$,
$E_{1/\lambda _j}$ isotropic for $\lambda _{j}\ne 1,-1$ and all the
$1,-1,\lambda _j,1/\lambda _j$ different. Let $A$ be an operator leaving
all the $E_{(..)}$ invariant, with ${A_\vert}_{A_\lambda }=\lambda
+N_\lambda $, $\lambda =\pm 1,\lambda _j,1/\lambda _j$ and $N_\lambda $
nilpotent. Then $A$ will be symplectic if (\Sy{}.3--5) hold.

\par We next consider logarithms of a symplectic matrix $A$. Decompose
${\bf C}^{2n}$ into generalized eigen-spaces as in (\Sy{}.2). We construct
$\log A$ with the same generalized eigen-spaces in the following way: On
$E_1$ we have
$A=1+N_1$ with
$N_1$ nilpotent and we put
$$\log A=\log (1+N_1)=N_1-{1\over 2}N_1^2+{1\over 3}N_1^3+...,$$
where the sum is finite, and in the following we always define the log
of $1+N$ in this way, when $N$ is nilpotent. If $\lambda \in\{\lambda
_j,\lambda _j^{-1}, -1\}$, we choose $\mu =\mu (\lambda )$ with $\lambda
=e^{\mu }$ in such a way that
\ekv{\Sy{}.6}{\mu ({1\over \lambda _j})=-\mu (\lambda _j).}
Write
$${A_\vert}_{E_\lambda }=\lambda +N_\lambda =\lambda (1+{1\over \lambda
}N_\lambda ),$$
and define on $E_\lambda $:
$$\log A=\mu +\log (1+{1\over \lambda }N_\lambda ).$$
This gives a definition which only depends on a choice of logarithms of
$\lambda _j$, $1\le j\le k$, and on $\log (-1)$ if $E_{-1}$ has positive
dimension. It is easy to check that
\ekv{\Sy{}.7}
{\log A+{\strans\log A}=(2k+1)2\pi i\Pi _{-1},}
for some integer $k$, where $\Pi _{-1}$ denotes the spectral projection
onto $E_{-1}$. Of course we have that $\exp\log A=A$. Recall also
that if
$B+{\strans B}=0$, then $\exp B$ is a symplectic matrix.

\par Assume now in addition that $A$ is a real matrix: $A:{\bf
R}^{2n}\to {\bf R}^{2n}$. Then $E_{\pm 1}$ become real in the sense that
they are invariant under complex conjugation $\Gamma :(x,\xi
)\mapsto\overline{(x,\xi )}$. The same holds for $E_\lambda $ if
$\lambda $ is real. If $\lambda $ is not real, there are two
possibilities:
\smallskip
\par\noindent 1) $\vert \lambda \vert \ne 1$. Then $\overline{\lambda
},\,{1\over \lambda },\, {1\over \overline{\lambda }}$ are also \ev{}s
and
$$E_\lambda \oplus E_{1\over \lambda }\oplus E_{\overline{\lambda
}}\oplus E_{{1\over \overline{\lambda }}}$$
is the complexification of a real symplectic space.\smallskip
\par\noindent 2) $\vert \lambda \vert =1$. Then $1/\lambda
=\overline{\lambda }$ and $E_\lambda \oplus E_{1\over \lambda }$ is
the complexification of a real symplectic space.
\smallskip

\par In all cases we have
$${A_\vert}_{E_{\overline{\lambda }}}=\Gamma ({A_\vert }_{E_\lambda
})\Gamma ,$$ and it is easy to see that $\log A$ will enjoy the same
property, provided that we choose $\mu (\lambda )$ in such a way that
\ekv{\Sy{}.8}{\mu (\overline{\lambda })=\overline{\mu (\lambda )}.}
This is possible, if we assume that $A$ has no negative \ev{}s. We get
\medskip

\par\noindent \bf \Prop{} \Sy{}.2. \it a) Let $A$ be a complex symplectic
$2n$-matrix, and choose a value $\mu (\lambda )$ of the logarithm of each
distinct \ev{} $\lambda $, different from 1 in such a way that (\Sy{}.6)
holds. Then we have a corresponding choice of $B=\log A$ with $e^B=A$,
satisfying (\Sy{}.7).\smallskip

\par\noindent b) Assume in addition that $A$ is real and has no negative
\ev{}s. Then by choosing $\mu (\lambda )$ with the additional property
(\Sy{}.8), $B=\log A$ becomes a real matrix, and ${\strans
B}+B=0$.\rm\medskip

\par From now on, we work under the assumptions of b) above, so that
$B=\log A$ is real and symplectically anti-symmetric. Consider the
quadratic form
\ekv{\Sy{}.9}
{b(\rho )={1\over 2}\sigma (\rho ,B\rho ).}
Then
$$\langle db(\rho ),t\rangle ={1\over 2}(\sigma (t,B\rho )+\sigma (\rho
,Bt))=\sigma (t,B\rho ).$$
On the other hand
$$\langle db(\rho ),t\rangle =\sigma (t,H_b(\rho )),$$
where $H_b$ denotes the Hamilton \vf{} associated to the function $b$,
so
\ekv{\Sy{}.10}
{H_b(\rho )=B\rho ,}
and the fact that $B=\log A$ can be expressed by
\ekv{\Sy{}.11}
{
A=\exp H_b .
}

\par We now consider the \pb{} of finding the "logarithm" of a
\canform{} also in the non-linear case. We first proceed somewhat
formally and let $p_s$ be a smooth real \fu{} depending smoothly on the
real parameter $s$. Consider the corresponding \canform{}
\ekv{\Sy{}.12}
{\kappa _{t,s}=\exp tH_{p_s}.}
We will later assume that $p_s$ vanishes to second order at some point
$\rho _0$, and then the germ of $\kappa _{t,s}$ at $\rho _0$ will be
well-defined for all real $t$. We differentiate the identity
$$\partial _t\kappa _{t,s}(\rho )=H_{p_s}(\kappa _{t,s}(\rho )),$$
\wrt{} $s$:
\ekv{\Sy{}.13}
{
\partial _t(\partial _s \kappa _{t,s}(\rho ))-({\partial H_{p_s}\over
\partial \rho }(\kappa _{t,s}(\rho ))(\partial _s\kappa _{t,s}(\rho
))=(\partial _sH_{p_s})(\kappa _{t,s}(\rho ))=H_{\partial _sp_s}(\kappa
_{t,s}(\rho )). }
Notice that the differential $d\kappa _{t,s}(\rho )\nu ={\partial \kappa
_{t,s}\over \partial \rho }(\rho )\nu $ satisfies
\ekv{\Sy{}.14}
{\partial _td\kappa _{t,s}(\rho )-{\partial H_{p_s}\over \partial \rho
}(\kappa _{t,s}(\rho ))\circ d\kappa _{t,s}(\rho )=0,\ d\kappa
_{0,s}(\rho )=1.}
Comparing the last two identities, we see that
\ekv{\Sy{}.15}
{\partial _s\kappa _{t,s}(\rho)=\int_0^t d(\kappa
_{t-\tilde{t},s})(\kappa _{\tilde{t},s}(\rho ))(H_{\partial _sp_s})(\kappa
_{\tilde{t},s}(\rho ))d\widetilde{t},}
which can also be written as
\ekv{\Sy{}.16}
{
\partial _s\kappa _{t,s}=\int_0^t (\kappa
_{t-\tilde{t},s})_*(H_{\partial _sp_s})d\widetilde{t}, }
where we use standard notation: lower $*$ for push forward and upper $*$
for pull back.

\par Rewrite (\Sy{}.16) as
$$\partial _s\kappa _{t,s}=(\kappa _{t,s})_*\int_0^t (\kappa
_{-\tilde t,s})_* H_{\partial _sp_s}d\widetilde{t},$$
and notice that
$$(\kappa _{-\tilde{t},s})_*H_{\partial _sp_s}=H_{(\kappa
_{-\tilde{t},s})_*\partial _sp_s}=H_{\partial _sp_s\circ \kappa
_{\tilde{t},s}},$$
since $\kappa _{-\tilde{t},s}$ is a \canform{}. Then
$$\partial _s\kappa _{t,s}=(\kappa _{t,s})_*\int_0^t H_{\partial
_sp_s\circ \kappa _{\tilde{t},s}}d\widetilde{t}=(\kappa
_{t,s})_*H_{\int_0^t \partial _sp_s\circ \kappa
_{\tilde{t},s}d\tilde{t}},$$
so
\ekv{\Sy{}.17}
{\partial _s\kappa _{t,s}=(\kappa _{t,s})_*H_{q_{t,s}},}
where
\ekv{\Sy{}.18}
{q_{t,s}=\int_0^t \partial _sp_s\circ \kappa
_{\tilde{t},s}d\widetilde{t}.}

\par In the last formula we shall take $t=1$ and consider a problem
where $\partial _sp_s$ will be the unknown. More precisely, let $\kappa
$ be a smooth \canform{}: ${\rm neigh\,}(0,{\bf R}^{2n})\to {\rm
neigh\,}(0,{\bf R}^{2n})$ with $\kappa (0)=0$. Let $A:=d\kappa
(0)=:\kappa _0$ have no negative \ev{}s so that part b) of \Prop{} \Sy{}.2
 applies. (More assumptions will be added later.) Let $B$ be a real
logarithm of $A$ as in the \prop{} and define the quadratic form
$p_0=b$ as in (\Sy{}.9). Then
\ekv{\Sy{}.19}
{\kappa _0=\exp H_{p_0}.}
We look for $p(\rho )=p_0(\rho )+{\cal O}(\rho ^3)$, so that
\ekv{\Sy{}.20}
{\kappa (\rho )=\exp H_p(\rho )+{\cal O}(\rho ^\infty ).
}
Let $\kappa _s$, $0\le s\le 1$, be a smooth family of \canform{}s with
\ekv{\Sy{}.21}
{
\kappa _s(0)=0
,\,\, d\kappa _s(0)=d\kappa (0),\,\, \kappa _0=d\kappa (0),\,\,\kappa
_1=\kappa .}
Then we look for a corresponding smooth family $p_s(\rho )=p_0(\rho
)+{\cal O}(\rho ^3)$, with
$p_{s=0}=p_0$ as above, such that
\ekv{\Sy{}.22}
{
\kappa _s(\rho )=\exp H_{p_s}(\rho )+{\cal O}(\rho ^\infty ).
}
Then $p=p_1$ will be a solution to our \pb{}. Define $q_s(\rho )={\cal
O}(\rho ^3)$, by $\partial _s\kappa _s=(\kappa _s)_*H_{q_s}$, or:
\ekv{\Sy{}.23}
{
(\kappa _s)^*\partial _s\kappa _s=H_{q_s}.
}
The discussion leading to (\Sy{}.18) indicates that we should find $p_s$
with the above properties, so that
\ekv{\Sy{}.24}
{
q_s(\rho )=\int_0^1 \partial _sp_s\circ \exp tH_{p_s}(\rho )dt+{\cal
O}(\rho ^\infty ). }

\par Let $N\ge 2$ and suppose that we have already found a smooth family
$p_s^{(N)}(\rho )=p_0(\rho )+{\cal O}(\rho ^3)$ with $p_0^{(N)}=p_0$, so
that
\ekv{{\rm E}_N}
{q_s(\rho )=\int_0^1 \partial _sp_s^{(N)}\circ \exp tH_{p_s^{(N)}}(\rho
)dt+{\cal O}(\rho ^{N+1}).}
Notice that $p_s^{(2)}=p_0$ solves (${\rm E}_2$) since $q_s(\rho )={\cal
O}(\rho ^3)$. Look for
$p_s^{(N+1)}=p_s^{(N)}+r_s^{(N+1)}$, with $r_s^{(N+1)}(\rho )={\cal
O}(s\rho ^{N+1})$. Then
$$\eqalign{\exp tH_{p_s^{(N+1)}}(\rho )&=\exp tH_{p_s^{(N)}}(\rho
)+{\cal O}(\rho ^N),\cr
\exp tH_{p_s^{(N+1)}}(\rho )&=\exp tH_{p_0}(\rho
)+{\cal O}(\rho ^2),
}$$
and we get
$$\eqalign{
&\int_0^1 \partial _sp_s^{(N+1)}\circ \exp tH_{p_s^{(N+1)}}(\rho)dt=\cr
&\int_0^1 \partial _sp_s^{(N)}\circ \exp tH_{p_s^{(N)}}(\rho)dt+
\int_0^1 \partial _sr_s^{(N+1)}\circ \exp tH_{p_0}(\rho)dt+{\cal
O}(\rho ^{N+2}). }$$
If we write the remainder in (${\rm E}_N)$ as $-v^{(N+1)}(\rho )+{\cal
O}(\rho ^{N+2})$, where $v^{(N+1)}$ is a homogeneous \pol{} of degree
$N+1$ depending smoothly on $s$, we will get (${\rm E}_{N+1}$), if we can
find
$\partial _sr_s^{(N+1)}$ as a homogeneous \pol{} $u^{(N+1)}$ of degree
$N+1$ depending smoothly on $s$, such that
\ekv{\Sy{}.25}
{
\int_0^1 u^{(N+1)}\circ \exp tH_{p_0}(\rho ) dt=v^{(N+1)}(\rho ).
}

\par Consider a general linear map $B:{\bf C}^{2n}\to {\bf C}^{2n}$ with
Jordan decomposition $B=D+N$, where $D$ is diagonalizable, $={\rm
diag\,}(d_j)$ \wrt{} a suitable basis, and $N$ is nilpotent and
commutes with $D$. The action of the \vf{} $Bx\cdot \partial _x$ on the
space $({\bf C}^{2n})^*$ of linear forms on ${\bf C}^{2n}$ can then be
identified with $\trans B$ in the natural way. Notice that $\trans B$
has the Jordan decomposition $\trans D +\trans N$ and that $\trans D$
becomes ${\rm diag\,}(d_j)$ if we select the dual basis to the one where
$D$ is diagonal. Define
$$\trans B^{(m)}={\trans B}\otimes 1\otimes .. \otimes 1+1\otimes {\trans
B}\otimes 1\otimes ..\otimes 1+...1\otimes .. \otimes 1\otimes{\trans
B},$$ as a linear endomorphism of the $m$-fold tensor product $(({\bf
C}^{2n})^*)^{\otimes m}$. We have the Jordan decomposition
\ekv{\Sy{}.26}
{\trans B^{(m)}={\trans D}^{(m)}+{\trans N}^{(m)}}
(where the first term to the right is diagonalizable, the second nilpotent
and the two terms commute). The corresponding \ev{}s are
$d_{j_1}+..+d_{j_m}$, for $j=(j_1,..,j_m)\in \{ 1,2,..,2n\}^{\{
1,2,..,m\} }$. The three operators in (\Sy{}.26) act naturally on the
symmetric tensor product $(({\bf C}^{2n})^*)^{\odot m}$ and the
decomposition (\Sy{}.26) is still a Jordan one on that space. The \ev{}s of
$\trans B^{(m)}$ become
\ekv{\Sy{}.27}
{
\beta _k=k_1d_1+..+k_{2n}d_{2n},\ k\in{\bf N}^{2n},\,\, k_1+..+k_{2n}=m.
}
$(({\bf C}^{2n})^*)^{\odot m}$ is equal to the space ${\cal P}^m_{{\rm
hom}}({\bf C}^{2n})$ of $m$-homogeneous polynomials on ${\bf C}^{2n}$
and $\trans B^{(m)}$ is the action of $Bx\cdot \partial _x$ on that
space.

\par Consider the map
\ekv{\Sy{}.28}
{
{\cal P}_{{\rm hom}}^m({\bf C}^{2n})\ni u\mapsto \int_0^1 u\circ \exp
(tB) dt\in {\cal P}_{{\rm hom}}^m({\bf C}^{2n}), }
which is equal to
\ekv{\Sy{}.29}
{
\int_0^1\exp (t{\trans B^{(m)}})dt.
}
The Jordan decomposition (\Sy{}.26) gives a similar decomposition of
(\Sy{}.29). The \ev{}s of (\Sy{}.29) are therefore
given by
\ekv{\Sy{}.30}
{
\int_0^1 e^{t\beta _k}dt=\cases{\displaystyle 1\hbox{ for }\beta _k=0\cr
\displaystyle {e^{\beta _k}-1\over \beta _k}\hbox{ for }\beta _k\ne 0},\
k\in{\bf N}^{2n},\,\, \vert k\vert =m.}
We conclude that the map (\Sy{}.28)
is invertible for a given $m$ precisely when for all $k\in{\bf N}^{2n}$
with $\vert k\vert =m$:
\ekv{\Sy{}.31}
{k_1d_1+..+k_{2n}d_{2n}\in 2\pi i{\bf Z}\Rightarrow
k_1d_1+..+k_{2n}d_{2n}=0.}

\par Now return to the equation (\Sy{}.25), where $p_0(\rho )=b(\rho
)={1\over 2}\sigma (\rho ,B\rho )$ and $B$ is the logarithm of the real
symplectic matrix $A=d\kappa (0)$, obtained under the assumptions of
Proposition \Sy{}.2, part b). The \ev{}s of $B$ are then of the form 0
with some possibly vanishing multiplicity and $\mu _j,-\mu _j\ne 0$ with
equal multiplicity  $>0$. Here we arrange so that all the \ev{}s are
distinct for instance by taking $\mu _j$ with either $\Re \mu _j>0$ or
with $\Re \mu _j=0$ and $\Im \mu _j>0$. We also recall that our set
of \ev{}s is closed under complex conjugation. The assumption that
(\Sy{}.31) holds for all $m$, then amounts to the assumption that
\ekv{\Sy{}.32}
{
\sum k_j\mu _j\in 2\pi i{\bf Z}\Rightarrow \sum k_j\mu _j=0,
}
for all $k_1,..,k_r\in{\bf Z}$. Here $r \le n$ is the number of
distinct $\mu _j$. (We could also have chosen to repeat the \ev{}s
according to their multiplicity without changing (\Sy{}.32).)

\par We have practically finished the proof of the following version of
a theorem of Lewis--Sternberg:
\medskip
\par\noindent \bf \Th{} \Sy{}.3. \it Let $\kappa :{\rm neigh\,}(0,{\bf
R}^{2n})\to {\rm neigh\,}(0,{\bf R}^{2n})$ be a smooth \canform{}.
Assume that $d\kappa (0)$ has no negative \ev{}s. Let the distinct
\ev{}s of $d\kappa (0)$ be 1 (possibly with multiplicity 0) and $\lambda
_j$, $\lambda _j^{-1}$, $1\le j\le r$ with $\vert \lambda _j\vert >1$ or
with $\vert \lambda _j\vert =1$ and $0<{\rm arg\,}\lambda _j<\pi $.
Choose $\mu _j$ with $\lambda _j=e^{\mu _j}$, in such a way that
$\overline{\lambda }_j$ corresponds to $\overline{\mu }_j$ and let
$B=\log A$ be given as in part b of \Prop{} \Sy{}.2. Assume that (\Sy{}.32)
holds, and let $p_0(\rho )=b(\rho )$ be given in (\Sy{}.9).

\par Then there exists $p(\rho )\in C^\infty ({\rm neigh\,}(0,{\bf
R}^{2n});{\bf R})$ such that $p(\rho )=p_0(\rho )+{\cal O}(\rho ^3)$ and
\ekv{\Sy{}.33}
{
\kappa (\rho )=\exp H_p(\rho )+{\cal O}(\rho ^\infty ).
}
$p$ is uniquely determined by these properties (for a given choice of
$p_0$).\rm
\medskip
\par This result (at least the existence part) is extremely close to a
corresponding one for complex \canform{}s, due to Lewis--Sternberg
([St], \Th{} 1 and Corollary 1.1) and clearly stated in [Fr],
Th{\'e}or{\`e}me V.1.
\par\noindent \bf End of the proof. \rm We establish the existence of
$p$. Let $\kappa _s$, $0\le s\le 1$ be a smooth family of \canform{}s
with $\kappa _1=\kappa $, $\kappa _s(0)=0$, $d\kappa _s(0)=d\kappa (0)$
and with $\kappa _0$ linear ($=d\kappa (0)$). Define
$q_s$ by (\Sy{}.23). The preceding discussion gives us a smooth family
$p_s(\rho )\in C^\infty ({\rm neigh\,}(0,{\bf R}^{2n});{\bf R})$ with
$p_{s=0}=p_0$, such that if $\widetilde{\kappa }_s=\exp H_{p_s}$, then
\ekv{\Sy{}.34}
{\widetilde{\kappa }_s^*\partial _s\widetilde{\kappa }_s=H_{q_s}+{\cal
O}(\rho ^\infty ).}
A simple computation shows that (\Sy{}.23,34) can be written as
\ekv{\Sy{}.35}
{\partial _s\kappa _s^{-1}(\rho )=-H_{q_s}(\kappa _s^{-1}(\rho )),\
\partial _s\widetilde{\kappa }_s^{-1}(\rho )=-H_{q_s}(\widetilde{\kappa
}_s^{-1}(\rho ))+{\cal O}(\rho ^\infty ),}
and we also know that
$\kappa _0^{-1}=\widetilde{\kappa }_0^{-1}$. It follows that
$$\kappa _s^{-1}(\rho )=\widetilde{\kappa }_s^{-1}+{\cal O}(\rho ^\infty
),$$
and hence that
$$\kappa _s(\rho )=\widetilde{\kappa }_s(\rho )+{\cal O}(\rho ^\infty
).$$
Taking $s=1$ gives (\Sy{}.33) with $p=p_1$.

\par We next prove the uniqueness of the Taylor expansion of $p$. Let
$\widetilde{p}$ have the same properties as $p$, so that
$$\exp H_p(\rho )=\exp H_{\widetilde{p}}(\rho )+{\cal O}(\rho ^\infty
),\ \widetilde{p}(\rho )=p(\rho )+{\cal O}(\rho ^3).$$
Assume that $\widetilde{p}-p$ does not vanish to infinite order, so that
$\widetilde{p}=p+r$, where $r(\rho )={\cal O}(\rho ^m)$, $r(\rho )\ne
{\cal O}(\rho ^{m+1})$, for some $3\le m\in{\bf N}$.

\par Put $p_s=(1-s)p+s\widetilde{p}=p+sr$, $0\le s\le 1$, so that
$p_0=p$, $p_1=\widetilde{p}$, and define $\kappa _s=\exp H_{p_s}$, so
that
\ekv{\Sy{}.36}
{
\kappa _1(\rho )=\kappa _0(\rho )+{\cal O}(\rho ^\infty ).
}
For this family, define $q_s$ by (\Sy{}.23). Then (\Sy{}.24) holds and since
$\partial _sp_s=r$, we have
$$q_s=\int_0^1 r\circ \exp (tH_{p_0})dt+{\cal O}(\rho ^{m+1}).$$
The previous discussion shows that the integral has a non-zero
Taylor polynomial of degree $m$:
\ekv{\Sy{}.37}
{
q_s(\rho )=\widetilde{q}(\rho )+{\cal O}(\rho ^{m+1}),\ 0\ne
\widetilde{q}\in{\cal P}^m_{{\rm hom}}. }

\par From (\Sy{}.35), we conclude that
$$\kappa _1^{-1}(\rho )-\kappa _0^{-1}(\rho )=-H_{\widetilde{q}}(\kappa
_0^{-1}(\rho ))+{\cal O}(\rho ^m),$$
which contradicts (\Sy{}.36). The proof is complete.\hfill{$\#$}\medskip

\bigskip

\centerline{\bf \Eq{}. Notions of equivalence.} 
\medskip
\par As in [Ia], our results will be valid "to infinite order at (0,0)"
and in this section we review the corresponding notions of equivalence.
Using these notions we also develop a very rudimentary functional
calculus for functions of several \pop{}s.

\par If $V_j\subset {\bf R}^n$ are open \neigh{}s of $0$ and $v_j\in
C^\infty (V_j)$, we say that $v_1$ and $v_2$ are equivalent; $v_1\equiv
v_2$, if $v_1-v_2$ vanishes to infinite order at $0$:
$v_1(x)-v_2(x)={\cal O}(x^\infty )$. This is clearly an equivalence
relation and the equivalence classes can be identified with the
corresponding formal Taylor expansions.

\par With $V=V_j\subset {\bf R}^n$ as above, let $S_{{\rm cl}}^0(V)$
denote the space of functions $a(x;h)$ in $C^\infty (V)$ depending on
the semi-classical parameter $h\in ]0,h_0]$  for some $h_0>0$, such that
\ekv{\Eq{}.1}{a(x,h)\sim \sum_0^\infty  h^ja_j(x),\ h\to 0,}
for some sequence of $a_j\in C^\infty (V)$. We say that $a^{(k)}\in
S_{{\rm cl}}^0(V_k)$, $k=1,2$ are equivalent and write $a^{(1)}\equiv
a^{(2)}$, if $a_j^{(1)}\equiv a_j^{(2)}$ for the corresponding
coefficients in (\Eq{}.1). Equivalently we can say that $a^{(k)}$ are
equivalent if $a^{(1)}(x;h)-a^{(2)}(x;h)={\cal O}((x,h)^\infty )$.

\par If $m>0$ is a smooth weight function on $V$, we define
$S^0(V,m)$ to be the space of smooth \fu{}s $a$ on $V$ such that for all
multi-indices $\alpha$, we have $\vert \partial _x^\alpha a(x)\vert
\le C_{\alpha ,a }m(x)$. Let $S^0_{{\rm cl}}(V,m)$ be the subspace of
$S^0_{{\rm cl}}$ for which (\Eq{}.1) holds in $S^0(V,m)$.

\par We next pass to the case of \pop{}s. Recall that if $p(x,\xi )$
belongs to an appropriate symbol class of functions on ${\bf R}^{2n}$,
then we define the corresponding $h$-Weyl quantization $P=p^w(x,hD_x)$
by:
\ekv{\Eq{}.2}
{
Pu(x)={1\over (2\pi )^n}\iint e^{i(x-y)\cdot \theta /h}p({x+y\over
2},\theta )u(y) dyd\theta . }
Recall that $p$ is called the Weyl-symbol of $P$. (See for instance
[DiSj].) Let
$S^0({\bf R}^{2n})$ denote the space of smooth functions that are \bdd{}
together with all their derivatives. If $p^{(k)}\in (S^0_{{\rm cl}}\cap
S^0)({\bf R}^{2n})=S_{{\rm cl}}^0({\bf R}^{2n},1)$, $k=1,2$, we say that
$(p^{(k)})^w(x,hD;h)$
 are equivalent (and use the symbol $\equiv$) if $p^{(k)}$
are equivalent in the sense of the classes $S^0_{{\rm cl}}$.

\par We will use the abbreviation ${\rm neigh\,}(0,{\bf R}^n)$ to denote
some \neigh{} of $0$ in ${\bf R}^n$. Tacitly it is understood that these
and other geometrical objects are \indep{} of $h$. We say that two
 smooth \canform{}s $\kappa _j:{\rm neigh\,}(0,{\bf R}^{2n})\to {\rm
neigh\,}(0,{\bf R}^{2n})$ with $\kappa _j(0)=0$ are equivalent if $\kappa
_1(\rho )-\kappa _2(\rho )={\cal O}(\rho ^\infty )$. Possibly after
shrinking the
\neigh{}s, we can introduce the inverses $\kappa _j^{-1}:{\rm
neigh\,}(0,{\bf R}^{2n})\to {\rm neigh\,}(0,{\bf R}^{2n})$. Then $\kappa
_1\equiv \kappa _{2}$ iff $\kappa _1^{-1}\equiv \kappa _2^{-1}$. Also
notice that the notion of equivalence of \canform{}s is stable under
composition in the natural way.

\par Let $\kappa :{\rm neigh\,}(0,{\bf R}^{2n})\to {\rm neigh\,}(0,{\bf
R}^{2n})$ be a \canform{} which maps $0$ to $0$. Then there exist
$N\in{\bf N}$ and a
\nondeg{} phase function $\phi (x,y,\theta )\in {\rm neigh\,}(0,{\bf
R}^{n+n+N})$ such that the graph of $\kappa $ in a \neigh{} of $(0,0)$
coincides with the image of the local \diff{}:
\ekv{\Eq{}.3}
{C_\phi :=\{ (x,y,\theta ); \phi '_\theta (x,y,\theta )=0\}\ni
(x,y,\theta )\mapsto (x,\phi _x';y,-\phi _y').}
Here we recall that a smooth real-valued function is called a
non-degenerate phase function (in the sense of H{\"o}rmander) if $d\phi
'_{\theta _1},..,d\phi _{\theta _N}'$ are linearly \indep{} on the set
$C_\phi $ above, which then becomes a $2n$-dimensional smooth
sub-\mfld{}. When discussing the relation between phases and symbols with
\canform{}s, it is tacitly understood that the point
$x=0,y=0,\theta =0$ corresponds to $\kappa (0)=0$ under the map (\Eq{}.3).

\par Let $\kappa $ be as above and let $\phi $ be a corresponding
generating phase. Let $\widetilde{\kappa }:\break {\rm neigh\,}(0,{\bf
R}^{2n})\to {\rm neigh\,}(0,{\bf R}^{2n})$ be a second \canform{} (with
the tacit convention that it also maps $0$ to $0$). It is easy to see
that $\widetilde{\kappa }\equiv \kappa $ if and only if
$\widetilde{\kappa }$ has a generating phase $\widetilde{\phi }$ which
is equivalent to $\phi $.

\par With $\phi ,\kappa $ as above we consider a \fop{} of order $0$:
\ekv{\Eq{}.4}
{Uu(x)=I(a,\phi )u(x)=h^{-{n+N\over 2}}\iint e^{{i\over h}\phi (x,y,\theta
)}a(x,y,\theta ;h)u(y)dyd\theta ,}
where $a\in S_{{\rm cl}}^0$ has its support in a \sufly{} small
\neigh{} of $(0,0,0)$. In this paper we only consider \fop{}s that are
elliptic at $(0,0,0)$  in the sense that $a_0(0,0,0)\ne 0$. In order to
normalize things, we will always assume that $\phi (0,0,0)=0$. If
$\kappa $ is the \canform{} generated by $\phi $, we say that $\kappa $
is the \canform{} associated to $U$. Thanks to the ellipticity assumption,
$\kappa $ is uniquely determined by $U$
in some \neigh{} of 0. We also recall the fundamental theorem about
\fop{}s, namely that if $\psi (x,y,w)$
 is a second phase which generates $\kappa $ and if $a$ has support in a
\sufly{} small \neigh{} of $(0,0,0)$, then there exists a classical
symbol $b(x,y,w;h)$ of order 0 with support in a small \neigh{} of
(0,0,0), such that $I(b,\psi )$ (formed as in (\Eq{}.4) with $N$
replaced by the dimension of $w$-space) is equal to $I(a,\phi
)$. Let $\widetilde{\kappa }$ be a second \canform{} with
$\widetilde{\kappa }\equiv \kappa $. Let $\widetilde{\phi }\equiv \phi $
be a corresponding generating phase. We say that
$\widetilde{U}=I(\widetilde{a},\widetilde{\phi })$ is equivalent to $U$
and write $\widetilde{U}\equiv U$, if $\widetilde{a}\equiv a$. It is a
standard exercise in \fop{} theory to verify that this definition of
equivalence does not depend on the choice $\phi $. It is also easy to
show that the definition is stable under composition in the natural way.

\par Below we will also need some functional calculus. First we consider
exponentials of \pop{}s. Let $p\sim \sum_0^\infty  p_j(x,\xi )h^j$ in
$S^0({\bf R}^{2n},1)$, and assume that $p_0$ is real-valued with
$p_0(0,0)=0$, $p_0'(0,0)=0$. Then
$e^{-itp^w(x,hD;h)/h}$ is well-defined for all complex $t$ (even without
the reality assumption on $p_0$) and for real $t$ we get a \fop{}. If
$\chi\in C_0^\infty ({\bf R}^{2n})$ is equal to 1 near $0$, then up to an
operator whose distribution kernel is rapidly decreasing together with
all its derivatives, we have that $\chi ^w e^{-itP/h}\chi ^w$ is a
\fop{} as above, with the associated \canform{} $\kappa _t=\exp tH_p$,
whose equivalence class does not depend on the choice of $\chi $. It is
also easy to see that if $\widetilde{P}$ is a second
\pop{} which is equivalent to $P$ and with real leading symbol, then for
real $t$, we have $e^{-itP/h}\equiv e^{-it\widetilde{P}/h}$ (in the
sense that we have equivalence for the corresponding truncated \op{}s).

\par Finally we discuss a very primitive pseudodifferential
functional calculus. Let $P_k=p_k(x,hD;h)$, $k=1,..,N_0$ be a commuting
\fy{} of \pop{}s with $p_k\in S_{{\rm cl}}^0({\bf R}^{2n},\langle (x,\xi
)\rangle ^m)$ (with the standard notation $\langle (x,\xi )\rangle
=(1+\vert (x,\xi )\vert ^2)^{1/2}$) and assume that the leading symbols
$p_{k,0}$ vanish at
$(0,0)$. Let $F(\iota_1,..,\iota_{N_0};h)\in S_{{\rm cl}}^0({\rm
neigh\,}(0,{\bf R}^{N_0}))$. Let $F_N$ be the sequence of polynomials in
$\iota_1,..,\iota_{N_0},h$ obtained by taking the Taylor polynomials of
order $N$ of the first $N$ terms in the \asy{} expansion of $F$, so that
$$F-F_N={\cal O}((\iota ,h)^N),\ (\iota,h)\to 0.$$
Then it is easy to see that $F_N(P_1,..,P_N;h)=q_N^w(x,hD_x;h)$, where
$$q_N(x,\xi ;h)-q_M(x,\xi ;h)={\cal O}((x,\xi ,h)^{k(N,M)}),\hbox{ where
}
k(N,M)\to \infty ,\, N,M\to \infty ,$$ and that this sequence defines
naturally an equivalence class of \pop{}s that we shall denote by
$F(P_1,..,P_{N_0};h)$.
\bigskip

\centerline{\bf \LF{}. Logarithms of \fop{}s.}
\medskip
\par Let $U_s$, $0\le s\le 1$ be a smooth family of elliptic \fop{}s of
order 0, associated to a fixed \canform{} $\kappa :{\rm neigh\,} (0,{\bf
R}^{2n})\to {\rm neigh\,}(0,{\bf R}^{2n})$ with $\kappa (0)=0$. We
represent $U_s$ by
\ekv{\LF{}.1}
{U_su(x)=h^{-{n+N\over 2}}\iint e^{{i\over h}\phi (x,y,\theta
)}u_s(x,y,\theta ;h)u(y)dyd\theta ,}
where $u_s\in S^0_{{\rm cl}}$ and more generally $\partial _s^ku_s\in
S^0_{{\rm cl}}$ for all
$k\in{\bf N}$  is a smooth family of classical symbols of order 0,
defined  in ${\rm neigh\,}((0,0,0);{\bf R}^{2n+N})$ and $\phi
$ is a real phase function which is \nondeg{} in the sense of H{\"o}rmander
[H{\"o}] and generates
$\kappa
$, so that $C_\phi \ni (x,y,\theta )\mapsto (x,\phi _x';y,-\phi
_y')\in{\rm graph\,}\kappa $ is a local \diffeo{}, where $C_\phi \subset
{\bf R}^{2n+N}$ is the sub-\mfld{} given by $\phi _\theta '(x,y,\theta
)=0$. To normalize things, we assume that $\phi '(0,0,0)=0$ and that
\ekv{\LF{}.2}{\phi (0,0,0)=0.}
Notice that this last assumption does not depend on the choice of phase
in the representation (\LF{}.1). In the following, we shall use the
equivalence relations "$\equiv$", defined in section\Eq{}.

\par We define the "logarithmic derivative" of our family, to be the
smooth family of \pop{}s $Q_s$ given by
\ekv{\LF{}.3}
{Q_s\equiv U_s^{-1}hD_sU_s.}
$Q_s$ and more generally $\partial _s^kQ_s$ is a smooth family of
classical \pop{}s defined in ${\rm neigh\,}((0,0);{\bf R}^{2n})$. (We
made an arbitrary choice of the order of the factors in (\LF{}.3), if
$\widetilde{Q}_s\equiv (hD_sU_s)U_s^{-1}$, then we get a new \pop{} which is
related to $Q_s$ by the intertwining relation
$U_sQ_s\equiv \widetilde{Q}_sU_s$.)

\par The family $U_s$ is determined uniquely by $U_0$ and its
logarithmic derivative:\medskip
\par\noindent \bf Lemma \LF{}.1. \it Let $V_s$ be a second family of
\fop{}s with the same properties as $U_s$ and associated to the same
\canform{} $\kappa $. Assume that $U_s^{-1}hD_sU_s\equiv V_s^{-1}hD_sV_s$ and
that $U_0\equiv V_0$. Then $U_s\equiv V_s$.\rm\medskip

\par\noindent \bf Proof. \rm Let $U$ be a fixed elliptic \fop{}
associated to $\kappa $, so that
$$U_s\equiv UA_s,\ V_s\equiv UB_s,$$
where $A_s$, $B_s$ are smooth families of \pop{}s. Then
$$U_s^{-1}hD_sU_s\equiv A_s^{-1}hD_sA_s,$$
and similarly for $V_s$, so we get
\ekv{\LF{}.4}
{
A_s^{-1}hD_sA_s\equiv B_s^{-1}hD_sB_s,\ A_0\equiv B_0.
}
>From this we conclude first that $A_s$ and $B_s$ have equivalent
principal symbols, then equivalent sub-principal symbols and so on, so
$A_s\equiv B_s$ and hence $U_s\equiv V_s$.\hfill{$\#$}\medskip
\par\noindent \it Remark. \rm We have $hD_sU_s\equiv U_sQ_s$, hence
$hD_sU_s^*\equiv -Q_s^*U_s^*$ for the adjoint \op{}s, so
$$hD_s(U_s^*U_s)+(Q_s^*(U_s^*U_s)-(U_s^*U_s)Q_s)\equiv 0.$$
$U_s^*U_s$ is a smooth family of elliptic \pop{}s and we conclude
\smallskip
\par\noindent a) If $U_s$ are unitary (up to equivalence), then $Q_s$ are
\sa{} (up to equivalence).\smallskip
\par\noindent b) If $U_s$ is unitary for one value of $s$ and $Q_s$ are
\sa{} for all $s$, then $U_s$ is unitary for all $s$ (again up
to equivalence).\medskip
 \par We now assume
for a while that
\ekv{\LF{}.5}
{U_s=U_{1,s},\hbox{ where }U_{t,s}\equiv e^{-itP_s/h},}
and $P_s$ is a smooth family of \pop{}s with the leading symbol $p(x,\xi
)$ \indep{} of $s$, so that $\kappa \equiv \exp H_p$ and $p(0,0)=0$ (thanks to
(\LF{}.2) and $p'(0,0)=0$ (since $\kappa (0,0)=(0,0)$). We shall derive a
simple formula for the logarithmic derivative: Start with
\ekv{\LF{}.6}
{hD_{\widetilde{t}}U_{\widetilde{t},s}+P_sU_{\widetilde{t},s}\equiv 0,\
U_{0,s}\equiv 1,}
and recall that $P_s$ and $U_{t,s}$ commute. Apply $hD_s$
to this relation:
$$hD_{\widetilde{t}}(hD_sU_{\widetilde{t},s})+P_s(hD_sU_{\widetilde{t},s})
\equiv -(hD_sP_s)U_{\widetilde{t},s},$$
which implies
$$hD_{\widetilde{t}}(U_{t-\widetilde{t},s}hD_s
U_{\widetilde{t},s})\equiv U_{t-\widetilde{t},s}(P_s+hD_{\widetilde{t}})
hD_sU_{\widetilde{t},s}\equiv  -U_{t-\widetilde{t},s}(hD_sP_s)
U_{\widetilde{t},s}.$$
Integrate this from $\widetilde{t}=0$ to $\widetilde{t}=t$:
$$hD_sU_{t,s}\equiv -{i\over h}\int_0^t U_{t-\widetilde{t},s}
(hD_sP_s)U_{\widetilde{t},s}d\widetilde{t}\equiv  -\int_0^t
U_{t-\widetilde{t},s}(\partial _sP_s)U_{\widetilde{t},s}d\widetilde{t}.$$
Taking $t=1$, we get the promised formula:
\ekv{\LF{}.7}
{
U_s^{-1}(hD_sU_s)\equiv -\int_0^1 U_{-t,s}(\partial _sP_s)U_{t,s}dt,
}
under the assumption (\LF{}.5).

\par Let $\kappa :{\rm neigh\,}(0,{\bf R}^n)\to {\rm neigh\,}(0,{\bf
R}^{2n})$ be a \canform{} as in \Th{} \Sy{}.3 (so that (\Sy{}.32) holds), and
choose $p=p_0+{\cal O}(\rho ^3)$ satisfying (\Sy{}.33):
\ekv{\LF{}.8}
{
\kappa (\rho )=\exp H_p(\rho )+{\cal O}(\rho ^\infty ).
}
Recall that $p$ is uniquely determined modulo ${\cal O}(\rho ^\infty )$
by $\kappa $ and the choice of the quadratic form $p_0$ with $\exp
H_{p_0}=d\kappa (0)$.

\par Let $U$ be an elliptic \fop{} of order 0 associated to the
\canform{} $\kappa $. We look for a \pop{} $P$ with leading symbol $p$
such that
\ekv{\LF{}.9}
{
U\equiv e^{-iP/h}
}

\par Let $P_0$ be a \pop{} with leading symbol $p$ and put
\ekv{\LF{}.10}
{U_0\equiv e^{-iP_0/h}.}
Let $[0,1]\ni s\mapsto U_s$ be a smooth family of \fop{}s as above, all
associated to $\kappa $ (modulo equivalence) and with $U_{s=0}=U_0$,
$U_1=U$. We look for a corresponding smooth \fy{} of \pop{}s $P_s$, with
leading symbol $p$, such that $P_{s=0}=P_0$, and
\ekv{\LF{}.11}
{U_s\equiv e^{-iP_s/h}.}
Then $P=P_1$ will solve (\LF{}.9).

\par Since the $U_s$ are associated to the same \canform{}, the
logarithmic derivative
\ekv{\LF{}.12}
{Q_s\equiv U_s^{-1}hD_sU_s, }
will be of order $-1$ (i.e. ${\cal O}(h^{+1})$) with Weyl symbol:
\ekv{\LF{}.13}
{
Q_s(\rho ;h)\sim hq_{s,1}(\rho )+h^2q_{s,2}(\rho )+...\,\, .
}
Motivated by (\LF{}.7) we shall first look for a smooth \fy{} $P_s$ with
leading symbol $p$ and $P_{s=0}=P_0$, such that
\ekv{\LF{}.14}
{
Q_s\equiv -\int_0^1 e^{itP_s/h}(\partial _sP_s)e^{-itP_s/h}dt.
}
Denoting the Weyl symbol of $P_s$ by the same letter,
\ekv{\LF{}.15}
{
P_s(\rho ;h)=p(\rho )+hp_{s,1}(\rho )+h^2p_{s,2}(\rho )+...\,\,,
}
we first see that $p_{s,1}$ should solve
\ekv{\LF{}.16}
{
q_{s,1}(\rho )=-\int_0^1(\partial _sp_{s,1})\circ \exp (tH_p)dt+{\cal
O}(\rho ^\infty ). }
As in the proof of \Th{} \Sy{}.3, we see that (\LF{}.16) has a unique solution
$\partial _sp_{s,1}$ (${\rm mod\,}{\cal O}(\rho ^\infty )$) and since
$p_{0,1}$ is given by the choice of $P_0$, we get a unique choice of
$p_{s,1}$.

\par Proceeding inductively, we assume that we have found $P_s^{(m)}$
with symbol
\ekv{\LF{}.17}
{
P_s^{(m)}(\rho ;h)=\sum_{j=0}^m h^jp_{s,j}(\rho )+\sum_{m+1}^\infty
h^jp_{0,j}(\rho ), }
where $p_{s,0}=p$ and $h^jp_{0,j}$ are the terms in the \asy{} expansion
of $P_0(\rho ;h)$, such that
\ekv{\LF{}.18}
{
-\int_0^1 e^{itP^{(m)}_s/h}(\partial _sP_s^{(m)})e^{-itP^{(m)}_s/h}
dt\equiv Q_s+h^{m+1}R_{m+1,s}, }
where $R_{m+1,s}$ is of order 0 with leading symbol $r_{m+1,s}$. We just
saw how to obtain this for $m=1$.

\par If $A$ is a \pop{} of order $0$, we see that
$e^{itP_s^{(m)}/h}Ae^{-itP_s^{(m)}/h}$ will change by an \op{} of order
$\le -(m+1)$ if we modify $P_s^{(m)}$ by an \op{} of order $\le -(m+1)$,
for instance by passing to $P_s^{(m+1)}$. It follows that
$$e^{itP_s^{(m+1)}/h}\partial
_sP_s^{(m+1)}e^{-itP_s^{(m+1)}/h}=e^{itP_s^{(m)}/h}\partial
_sP_s^{(m+1)}e^{-itP_s^{(m)}/h}+{\cal O}(h^{m+2}).$$
To get $P_s^{(m+1)}$ satisfying (\LF{}.18) with $m$ replaced by $m+1$, it
suffices to have
\ekv{\LF{}.19}
{
\int_0^1 e^{itP_s^{(m)}/h}\partial
_s(P_s^{(m+1)}-P_s^{(m)})e^{-itP_s^{(m)}/h}dt \equiv  h^{m+1}R_{m+1,s}+{\cal
O}(h^{m+2}) }
(with the same $R_{m+1,s}$ as in (\LF{}.18)),
which gives for the leading symbols
\ekv{\LF{}.20}
{
\int_0^1 (\partial _sp_{s,m+1})\circ \exp (tH_p) dt \equiv r_{m+1,s}.
}
Again this has a unique solution $\partial _sp_{s,m+1}$, and our
induction procedure can be continued and gives a solution $P_s$ to
(\LF{}.14).

\par Let $\widetilde{U}_s=e^{-itP_s/h}$. Then by construction
$\widetilde{U}_0=U_0$, $\widetilde{U}_s^{-1}hD_s\widetilde{U}_s\equiv
U_s^{-1}hD_sU_s$, and Lemma \LF{}.1
 implies that $\widetilde{U}_s\equiv U_s$ and in particular that
\ekv{\LF{}.21}{U\equiv e^{-iP/h},\ P=P_1.}
This gives the existence part of the following \medskip
\par\noindent \bf \Th{} \LF{}.2. \it Let $\kappa :{\rm neigh\,}(0,{\bf
R}^{2n})\to {\rm neigh\,}(0,{\bf R}^{2n})$ be a smooth \canform{} as in
\Th{} \Sy{}.3 and choose $\mu _j$, $p_0$ as there, so that (\Sy{}.32) holds.
Let $p\in C^\infty ({\rm neigh\,}(0,{\bf R}^{2n});{\bf R})$ be the
unique function ${\rm mod\,}{\cal O}(\rho ^\infty )$ of the form
$p=p_0+{\cal O}(\rho ^3)$ with $\kappa (\rho )=\exp H_p(\rho )+{\cal
O}(\rho ^\infty )$.

\par Let $U$ be an elliptic \fop{} of order 0 associated to $\kappa $.
Then there exists a \pop{} $P^w(x,hD_x;h)$ with symbol
\ekv{\LF{}.22}
{
P(\rho ;h)\sim p(\rho )+hp_1(\rho )+...,
}
such that
\ekv{\LF{}.23}
{
U\equiv e^{-iP/h}.
}
$P$ is uniquely determined modulo "$\equiv$" and up to an integer
multiple of $2\pi h$ by (\LF{}.23) and the choice of
$p_0$.\rm\medskip

\par It remains to prove the uniqueness modulo "$\equiv$". Let
$\widetilde{P}^w(x,hD_x;h)$ be another \op{} with the same properties;
\ekv{\LF{}.24}
{
\widetilde{P}(\rho ;h)\sim \widetilde{p}(\rho )+h\widetilde{p}_1(\rho
)+...,\ \widetilde{p}=p_0+{\cal O}(\rho ^3). }
Then we must have $\kappa (\rho )=\exp H_{\widetilde{p}}(\rho )+{\cal
O}(\rho ^\infty )$ and from the uniqueness part of \Th{} \Sy{}.3, we
conclude that $\widetilde{p}=p+{\cal O}(\rho ^\infty )$.

\par Put $P_s=(1-s)P+s\widetilde{P}$, $0\le s\le 1$, so that $P_0=P$,
$P_1=\widetilde{P}$ and define $U_s=e^{-iP_s/h}$. For this family,
define $Q_s$ by (\LF{}.12). If $\widetilde{P}\not\equiv P$, let $1\le m\le
\infty $ be the smallest integer with
\ekv{\LF{}.25}
{\widetilde{p}_m-p_m\ne {\cal O}(\rho ^\infty ).}
If $m=1$, we may also assume that $\widetilde{p}_1-p_1$ is not $\equiv$
to an integer multiple of $2\pi $. From (\LF{}.14), we see that
$$Q_s\sim \sum_1^\infty  h^jq_{s,j},$$
with $q_{s,j}(\rho )={\cal O}(\rho ^\infty )$ for $1\le j\le m-1$, and
with
\ekv{\LF{}.26}
{q_{m,s}=-\int_0^1 (\widetilde{p}_m-p_m)\circ \exp (tH_p)dt+{\cal O}(\rho
^\infty ).}
When $m=1$, $q_{m,s}$ is not $\equiv$ to an integer multiple of $2\pi $.
>From the invertibility of the map (\Sy{}.28), we conclude that
\ekv{\LF{}.26}
{
q_m:=q_{m,0}=q_{m,s}+{\cal O}(\rho ^\infty ),\ q_m\ne {\cal
O}(\rho ^\infty ). }

\par Let $W_s$, $0\le s\le 1$ be a smooth \fy{} of \fop{}s which solves
\ekv{\LF{}.27}
{Q_s\equiv W_s^{-1}hD_sW_s,\ W_0=1.}
If $q_{0,s}$ had been 0 rather than just ${\cal O}(\rho ^\infty )$, the
$W_s$ would have been \pop{}s, so in general they are equivalent to such
\op{}s:
\ekv{\LF{}.28}
{
W_s\equiv R_s^w(x,hD_x;h),\hbox{ where }R_s(\rho ;h)\sim\sum_{j=0}^\infty
h^jr_{j,s}(\rho ), }
satisfying $$r_{0,s}(\rho )^{-1}\partial _sr_{0,s}(\rho )=iq_1$$ when
$m=1$ and
$$r_{0,s}(\rho )=1,\,\, r_{j,s}=0\hbox{ for }1\le j\le m-2,\ \partial
_sr_{m-1,s}=iq_m,$$
when $m\ge 2$. In other words,
$$\eqalign{r_{0,s}(\rho )&=e^{isq_1(\rho )},\hbox{ when }m=1,\cr
R_s(\rho ;h)&=1+isq_mh^{m-1}+{\cal O}(h^m),\hbox{ when }m\ge 2.
}$$
In both cases, we have $R_1\not\equiv 1$, so
\ekv{\LF{}.29}
{W_1\not\equiv 1.}

\par If we put $\widetilde{U}_s=U_0W_s$, we see that
$\widetilde{U}_0=U_0$ and that
$$Q_s\equiv \widetilde{U}_s^{-1}hD_s\widetilde{U}_s\equiv U_s^{-1}hD_s
U_s.$$
By Lemma \LF{}.1 we conclude that $\widetilde{U}_s\equiv U_s$ and in
particular,
$$U_0W_1=\widetilde{U}_1\equiv U_1.$$
Since $W_1\not\equiv 1$, this contradicts the assumption that $U_1\equiv
U_0$, and the proof of \Th{} \LF{}.2 is complete.\hfill{$\#$}\medskip

\par\noindent \it Remark. \rm Up to equivalence we have that $U$ is
unitary iff $P$ is \sa{}:
$$U^*U\equiv 1\Leftrightarrow P^*\equiv P.$$
Indeed (\LF{}.23) gives
$$(U^*)^{-1}\equiv e^{-iP^*/h},$$
and it suffices to apply the uniqueness statement in \Th{} \LF{}.2.\medskip

\bigskip

\centerline{\bf \BN{}. Birkhoff normal forms.}
\medskip
\par To get a normal form for the \fop{} in \Th{} \LF{}.2, it suffices to
get the quantized "Birkhoff" normal form of the \op{} $P$. For
simplicity we shall make a non-resonance assumption, and simply recall
how this was done in [Sj] in a slightly less general setting (in the
spirit of works of Bellissard--Vittot, Graffi--Paul and others cited
there). The extension to the present case is however completely
immediate.

\par Let $P\sim p(\rho )+hp_1(\rho )+..$ be as in (\LF{}.22), with $p$
real, $p(0)=0$, $p'(0)=0$. Put $p_0(\rho )={1\over 2}\langle p''(0)\rho
,\rho \rangle $ and let $B$ be the corresponding fundamental matrix, so
that
\ekv{\BN{}.1}
{p_0(\rho )={1\over 2}\sigma (\rho ,B\rho ),\ {\strans B}=-B.}
Let $\mu _j,-\mu _j$ and possibly $0$ be the distinct \ev{}s of $B$. We
recall that \Th{} \LF{}.2 was obtained under the assumption (\Sy{}.32).
We add a non-resonance assumption, and for that purpose we
temporarily change the notation slightly and denote by $\mu _j,-\mu _j$,
$1\le j\le n$ all the \ev{}s of $B$, possibly repeated according to their
multiplicity. Assume
\ekv{\BN{}.2}
{
\sum_1^n k_j\mu _j=0,\,\, k_j\in{\bf Z}\Rightarrow k_1=..=k_n=0.
}
This implies that the $\mu _j$ are distinct and $\ne 0$, so $B$ has the
$2n$ distinct \ev{}s $\mu _j,-\mu _j$, $1\le j\le n$, which is in
agreement with the earlier notation with $r=n$.

\par Notice that (\Sy{}.32) and (\BN{}.2) combine into the single
condition
\ekv{\BN{}.3}
{\sum_1^n k_j\mu _j\in 2\pi i{\bf Z},\ k_j\in{\bf
Z}\Rightarrow k_1=..=k_n=0,}
which does not change if we modify the choice of the $\mu _j$ by some
multiples of $2\pi i$.

\par Let $e_1,..,e_n,f_1,..,f_n\in{\bf C}^{2n}$ be a basis of
eigen-vectors of $B$, associated to $\mu _1,..,\mu _n$, $-\mu _1,..,-\mu
_n$. Then $\sigma (e_j,e_k)=\sigma (f_j,f_k)=0$
and $\sigma (f_j,e_k)=0$ for $j\ne k$. We can arrange so that
$$\sigma (f_j,e_k)=\delta _{j,k},$$
and then we have a symplectic basis in ${\bf C}^{2n}$. The corresponding
coordinates $x_j,\xi _j$ given by
${\bf C}^{2n}\ni \rho =\sum_1^n(x_je_j+\xi _jf_j)$ will be symplectic,
and in these coordinates, we get
\ekv{\BN{}.4}
{
p_0(\rho )=\sum_1^n \mu _jx_j\xi _j,
}
with the Hamilton field
\ekv{\BN{}.5}
{H_{p_0}=\sum_1^n\mu _j(x_j\partial _{x_j}-\xi _j\partial _{\xi _j}).}
If we consider $H_{p_0}:{\cal P}_{\rm hom}^m\to {\cal P}_{\rm hom}^m$,
we see that the monomials $x^\alpha \xi ^\beta $, $\vert \alpha \vert
+\vert \beta \vert =m$ form a basis of eigen-vectors and
\ekv{\BN{}.6}
{
H_{p_0}(x^\alpha \xi ^\beta )=\mu \cdot (\alpha -\beta )x^\alpha \xi
^\beta , }
where $\mu =(\mu _1,..,\mu _n)$. The assumption (\BN{}.2) implies that $\mu
\cdot (\alpha -\beta )=0$ precisely when $\alpha =\beta $, so if we
let the set of resonant polynomials ${\cal R}_{{\rm hom}}^m\subset {\cal
P}_{\rm hom}^m$ be the space of linear combinations of all the $x^\alpha
\xi ^\alpha =(x_1\xi _1)^{\alpha _1}...(x_n\xi _n)^{\alpha
_n}$
 with $2\vert \alpha \vert =m$, we see that $H_{p_0}$ induces a
bijection from ${\cal P}_{\rm hom}^m/{\cal R}_{\rm hom}^m$ into itself.

\par We say that $u\in C^\infty ({\rm neigh\,}(0,{\bf R}^{2n}))$ is
resonant if its Taylor expansion at 0 only contains resonant polynomials.
Since
$p_0$ is real it is easy to see that the space of resonant smooth \fu{}s
is closed under complex conjugation. We also see that $u$ is resonant
iff $\exists f\in C^\infty ({\rm neigh\,}(0,{\bf C}^n))$ with
$\overline{\partial }f(\rho )={\cal O}(\rho ^\infty )$ such that
$$u(x)=f(x_1\xi _1,..,x_n\xi _n)+{\cal O}(\rho ^\infty ).$$
Considering Taylor expansions it is easy to get (cf [Sj]):
\medskip
\par\noindent \bf Lemma \BN{}.1 \it For every $v\in C^\infty ({\rm
neigh\,}(0,{\bf R}^{2n}))$, $\exists u\in C^\infty ({\rm neigh\,}(0,{\bf
R}^{2n}))$ unique up to a resonant \fu{}, such that
$$H_{p_0}u=v+r,$$
where $r$ is resonant. If $v={\cal O}(\rho ^m)$, we can find $u,r$ with
the same property.\rm\medskip

\par As for $H_p$ we only give the corresponding existence statement:
\medskip
\par\noindent \bf Lemma \BN{}.2. \it For every $v\in C^\infty ({\rm
neigh\,}(0,{\bf R}^{2n}))$, $\exists u\in C^\infty ({\rm neigh\,}(0,{\bf
R}^{2n}))$, such that
$$H_pu=v+r,$$
where $r$ is resonant. If $v={\cal O}(\rho ^m)$, we can choose $u,r$
with the same property.\rm\medskip

\par Notice that since $p$ is real, if $v$ is real, we can take $u,r$
real. The classical Birkhoff normal form is then given in
\medskip
\par\noindent\bf \Prop{} \BN{}.3. \it $\exists$ a smooth \canform{} $\kappa
:{\rm neigh\,}(0,{\bf R}^{2n})\to {\rm neigh\,}(0,{\bf R}^{2n})$,
such that $\kappa (\rho )=\rho +{\cal O}(\rho ^2)$, and
$$p\circ \kappa =p_0+r,$$
where $r$ is resonant and ${\cal O}(\rho ^3).$
\medskip
\par\noindent \bf Proof. \rm If $q\in C^\infty ({\rm neigh\,}(0,{\bf
R}^{2n});{\bf R})$, $q(\rho )={\cal O}(\rho ^{m+1})$, with $m\ge 2$,
then we see that $\exp H_q(\rho )=\rho +H_q(\rho )+{\cal O}(\rho
^{2m-1})$. Let first $q_3={\cal O}(\rho ^3)$ solve
$H_p(q_3)=(p-p_0)-r_3$, where $r_3(\rho )={\cal O}(\rho ^3)$ is
resonant. Let $\kappa _2(\rho )=\exp H_{q_3}(\rho )$. Then
$$p(\kappa _2(\rho ))=p(\rho +H_{q_3}(\rho )+{\cal O}(\rho ^3))=p_0(\rho
)+r_3(\rho )+{\cal O}(\rho ^4)=:\widetilde{p}(\rho )+r_3(\rho ).$$
(We used Lemma 1.1 in [Sj].)
Now repeat the argument with $p$ replaced by $\widetilde{p}$ and find
$\kappa _3=\exp H_{q_4}$ et.c. Finally, we choose $\kappa $ with $\kappa
(\rho )\sim \kappa _2\circ \kappa _3\circ \kappa _{(..)}\circ ...$. See
for instance [Sj] for more details.\hfill{$\#$}\medskip

\par We next review the quantum normal form of a \pop{}. Let
$P_0=p_0^w(x,hD_x)$. If $A$ is (equivalent to) a \pop{} with symbol
$a\sim a_0(\rho )+ha_1(\rho )+...$, we say that $A$ is resonant if every
$a_j$ is resonant. Since $a_j$ is resonant precisely when
$H_{p_0}a_j={\cal O}(\rho ^\infty )$ and $[P_0,A]$ has the symbol
${h\over i}H_{p_0}a$, we see that $A$ is resonant if  $[P_0,A]\equiv
0$. (Later we shall also recall that $A$ is resonant precisely when it
is equivalent to a function of the elementary action \op{}s.)

\par With $p$ as above, let $P=P^w$ be a \pop{} with leading symbol $p$,
so that $P(\rho ;h)\sim p(\rho )+hp_1(\rho )+...$. Let $\kappa $ be as
in \Prop{} \BN{}.3 and let $U$ be a corresponding elliptic \fop{} that we
choose to be microlocally unitary near $0$. Then
$U^{-1}PU\equiv \widetilde{P}$, where $\widetilde{P}$ has the leading symbol
$\widetilde{p}=p_0+r$ with $r={\cal O}(\rho ^3)$ resonant. We drop
the tilde and continue the reduction of "$P=\widetilde{P}$" by means of
conjugation with \pop{}s. We look for a \pop{} $Q=Q^w$ of order 0, such
that
\ekv{\BN{}.7}
{e^{iQ}Pe^{-iQ}=P_0+R,}
where $R$ is resonant. Here the \lhs{} can also be written
\ekv{\BN{}.8}
{e^{iQ}Pe^{-iQ}=e^{i{\rm ad}_Q}P=P+i{\rm ad}_QP+{(i{\rm ad}_Q)^2\over
2}P+...,}
where the sum is asymptotic in $h$, since ${\rm ad}_Q^kP$ is of order
$\le -k$. We look for $Q$ with symbol $q_0+hq_1+...$. The leading symbol
of $i{\rm ad}_QP$ is $hH_{q_0}p=-hH_pq_0$, so we first choose $q_0$ so
that
\ekv{\BN{}.9}
{H_pq_0=p_1+r_1,}
with $r_1$ resonant. Then the first two terms in the \asy{} expression of
the \op{} (\BN{}.8) become resonant. The choice of $q_1$ will influence the
$h^2$ term in the symbol of $e^{iQ}Pe^{-iQ}$ only via the term $i{\rm
ad}_QP$, and to make the $h^2$ term resonant, leads to a new equation of
the same type as (\BN{}.9). It is clear that this construction can be
iterated and we find $Q$ so that (\BN{}.7) holds with $R$ resonant.

\par If the original symbol $P$
is \sa{}, then the new "$P=\widetilde{P}=U^{-1}PU$" will also be \sa{} and
hence have a real-valued symbol. We can then find $Q$ in (\BN{}.7) \sa{},
because of the observation that if $A,B$ are \sa{}, then $i{\rm ad}_AB$
is \sa{}, so if $Q,P$ are \sa{}, then all terms of the last expression
in (\BN{}.8) have the same property. Consequently, in each step of the
computation, we will encounter an equation of the form
$H_pq_k=\widehat{p}_k+r_k$, with $\widehat{p}_k$ real-valued, and we
then choose the solution $q_k$ and the resonant remainder $r_k$ to be
real.  This means that $e^{-iQ}$ will be unitary. If $V=Ue^{-iQ}$, we
finally obtain for the original $P$, that
$$V^{-1}PV\equiv P_0+R,$$
where $R$ is resonant of order 0 with leading symbol $r={\cal O}(\rho
^3)$. Summing up we have\medskip

\par\noindent \bf \Th{} \BN{}.4. \it Let $p(\rho )=p_0(\rho )+{\cal
O}(\rho ^3)$ be real-valued with $p_0(\rho )={1\over 2}\sigma (\rho
,B\rho )$, where $B$ is symplectically anti-symmetric satisfying the
non-\res{} condition (\BN{}.2). Let $P$ be a \pop{} with leading symbol
$p$. Then there exists an elliptic \fop{} $V$ associated to the
\canform{} $\kappa $ in \Prop{} \BN{}.3, such that
\ekv{\BN{}.10}
{
V^{-1}PV\equiv P_0+R
}
where $R$ is a resonant \pop{} of order $\le 0$ with leading symbol
$={\cal O}(\rho ^3)$. If $P$ is \sa{}, we can choose $V$ to be unitary.
\rm\medskip

\par When applying this to $U$ and $P$ in \Th{} \LF{}.2, we notice that
\ekv{\BN{}.11}
{
V^{-1}UV\equiv e^{-iV^{-1}PV/h},
}
which can be viewed as a quantum normal form for our \fop{} $U$.

\par In the appendix to this section, we review that under the
non-\res{} assumption (\BN{}.2), there are real symplectic coordinates
$x_1,..,x_n,\xi _1,..,\xi _n$ such that
\eekv{\BN{}.12}
{p_0(\rho )=\sum_1^{n_{\rm hc}}(\alpha _j(x_{2j-1}\xi _{2j-1}+x_{2j}\xi
_{2j})-\beta _{j}(x_{2j-1}\xi _{2j}-x_{2j}\xi _{2j-1}))}
{\hskip 1cm +\sum_{2n_{\rm hc}+1}^{2n_{\rm hc}+n_{\rm hr}}\mu _jx_j\xi
_j+\sum_{2n_{\rm hc}+n_{\rm hr}+1}^n\nu _j{1\over 2}(x_j^2+\xi _j^2),}
where
$\nu _j\in{\bf R}$ are non-vanishing with distinct values of $\vert
\nu _j\vert $, $\mu _j>0$ are distinct, and $\alpha _j,\beta _j>0$  with
$\alpha _j+i\beta _j$ distinct. We have the corresponding resonant
"actions":
\eeekv{\BN{}.13}
{\cases{\iota _j=x_{2j-1}\xi _{2j-1}+x_{2j}\xi _{2j}\cr
\iota_{j+n_{\rm hc}}=x_{2j-1}\xi _{2j}-x_{2j}\xi _{2j-1}},\ 1\le j\le
n_{\rm hc}, }
{\iota _j=x_j\xi _j,\ \ \ \ \ \ 2n_{\rm hc}+1\le j\le 2n_{\rm hc}+n_{\rm
hr},} {\iota_j={1\over 2}(x_j^2+\xi _j^2),\ \ \  2n_{\rm hc}+n_{\rm
hr}+1\le j\le 2_{\rm hc}+n_{\rm hr}+n_e=n.}

\par A resonant \fu{} is one which can be written
$f(\iota_1,..,\iota_n)+{\cal O}(\rho ^\infty )$ for some smooth \fu{}
$f$, and using the simple functional calculus of section \Eq{}, we see
that a \pop{}
$R$ of order 0 is resonant iff
$R\equiv F(I_1,..,I_n;h)$, where $F(\iota ;h)$ is a classical symbol of
order $0$ and $I_j=\iota_j^w(x,hD_x;h)$ is the corresponding commuting
family of quantized actions. (We refer to [DiSj] and references there to
the original work of B. Helffer and D. Robert, for more elaborate
functional calculi.) Combining this with
\Th{}
\BN{}.4 and (\BN{}.11), we get
\ekv{\BN{}.14}
{
V^{-1}UV\equiv e^{-iF(I_1,..,I_n;h)/h,}}
where
\ekv{\BN{}.15}
{F(\iota ;h)\sim\sum_0^\infty  F_j(\iota )h^j,}
and

\ekv{\BN{}.16}
{F_0(\iota )=\sum_1^{n_{\rm hc}}(\alpha _j\iota_j-\beta
_{j}\iota_{n_{\rm hc}+j})+\sum_{2n_{\rm hc}+1}^{2n_{\rm
hc}+n_{\rm hr}}\mu _j\iota_j+\sum_{2n_{\rm hc}+n_{\rm hr}+1}^n\nu
_j\iota_j+{\cal O}(\iota^2),}

\bigskip

\centerline{\bf Appendix. Review of the real normal form for the
quadratic part.}
\medskip
\par Here we review some standard facts. See also [Ze2], [It].
Let $B:{\bf R}^{2n}\to {\bf R}^{2n}$ be the symplectically
anti-symmetric matrix of \Prop{} \Sy{}.2, case b. We make the non-\res{}
assumption (\BN{}.2), so that all the \ev{}s of $B$ are simple and $\ne
0$. Recall from section \Sy{} that they can be grouped into families of 2
or 4 according to the following 3 cases:\smallskip

\par\noindent \it Case 1. \rm $\mu>0$ is an \ev{}. Then $-\mu $ is also
an \ev. Let $e,f$ be corresponding real eigen-vectors with $\sigma
(f,e)=1$, spanning a real symplectic space of dimension 2. A point in
this subspace can be written $\rho =xe+\xi f$, so that $(x,\xi )$ become
symplectic coordinates, and we get
$p_0(\rho )=b(\rho )={1\over 2}\sigma (\rho ,B\rho )={1\over 2}\sigma
(xe+\xi f,\mu xe-\mu \xi f)=\mu x\xi $.
The corresponding resonant action is $x\xi $.
\smallskip
\par\noindent
\it Case 2. \rm $\mu $ is an \ev{} with $\Re \mu ,\Im \mu >0$. Then
$-\mu ,\overline{\mu },-\overline{\mu }$ are also \ev{}s, and we let
$e,f,\overline{e},\overline{f}$ be corresponding eigen-vectors. We have
\ekv{{\rm A}.1}{
\sigma (e,\overline{e})=0,\ \sigma (e,\overline{f})=0,\ \sigma
(f,\overline{f})=0,
}
 and if $e$ is fixed, we can choose $f$ so that
\ekv{{\rm A}.2}{\sigma (f,e)=1.}
$e,f$ and $\overline{e},\overline{f}$ span 2-dimensional complex
symplectic subspaces that are complex conjugate to each other, while
$e,\overline{e},f,\overline{f}$ span a 4-dimensional symplectic space
which is the complexification of a real symplectic space (of real
dimension 4).

\par Writing $\rho =ze+\zeta f+w\overline{e}+\omega \overline{f},$
we get
\ekv{{\rm A}.3}
{b(\rho )={1\over 2}\sigma (\rho ,B\rho )=\mu z\zeta +\overline{\mu
}w\omega .}
The resonant action terms are $z\zeta $ and $w\omega $.

\par To get the real canonical form, we write
$$e={1\over \sqrt{2}}(e_1+ie_2),\ f={1\over \sqrt{2}}(f_1-if_2),$$
with $e_j,f_j$ real. Using this in (A.1,2), we get
\ekv{{\rm A}.4}
{\sigma (e_j,e_k)=\sigma (f_j,f_k)=0,\ \sigma (f_j,e_k)=\delta _{j,k},}
so $e_1,e_2,f_1,f_2$ is a symplectic basis in the real symplectic space
mentioned above.

\par We also have the inverse relations
$$\eqalign{
&e_1={1\over \sqrt{2}}(e+\overline{e}),\ e_2={1\over
i\sqrt{2}}(e-\overline{e}),\cr
&f_1={1\over \sqrt{2}}(f+\overline{f}),\ f_2={i\over
\sqrt{2}}(f-\overline{f}).  }$$
Write
$$\rho =ze+\zeta f+w\overline{e}+\omega
\overline{f}=\sum_1^2x_je_j+\sum_1^2 \xi _jf_j,$$
so that $(x,\xi )$ are real symplectic coordinates on our symplectic
4-space. Then
$$\eqalign{z&={1\over \sqrt{2}}(x_1-ix_2),\ w={1\over
\sqrt{2}}(x_1+ix_2),\cr
\zeta &={1\over \sqrt{2}}(\xi _1+i\xi _2),\ \omega ={1\over \sqrt{2}}(\xi
_1-i\xi _2),}$$
and using this in (A.3), we get
\ekv{{\rm A}.5}
{
b(\rho )=\alpha (x_1\xi _1+x_2\xi _2)-\beta (x_1\xi _2-x_2\xi _1),
}
with $\mu =\alpha +i\beta $.

\par The resonant actions can also be written
$$\eqalign{
z\zeta &={1\over 2}((x_1\xi _1+x_2\xi _2)+i(x_1\xi _2-x_2\xi
_1))\cr
w\omega &={1\over 2}((x_1\xi _1+x_2\xi _2)-i(x_1\xi _2-x_2\xi
_1)).}$$
The (resonant) real-valued functions (on the real symplectic 4-space
above) which only depend on $z\zeta ,w\omega $ are precisely the
functions of
$x_1\xi _1+x_2\xi _2$, $x_1\xi _2-x_2\xi _1$ modulo ${\cal O}(\rho
^\infty )$. Notice that these two functions Poisson commute.\smallskip

\par\noindent \it Case 3. \rm $\mu \ne 0$ is an \ev{} with $\Re \mu =0$.
Then $\overline{\mu }=-\mu $ is also an \ev{}. If $e$ is an eigen-vector
corresponding to $\mu $, then $\overline{e}$ will be an eigen-vector
corresponding to $\overline{\mu }$ and $\sigma (e,\overline{e})\in i{\bf
R}\setminus \{ 0\}$. Possibly after permuting $\mu $ and $-\mu $ and
after normalization, we can assume that
$${1\over i}\sigma (e,\overline{e})=1.$$
$e,\overline{e}$ span a 2-dimensional symplectic subspace which is the
complexification of a corresponding real 2-dimensional space. Let
$f=i\overline{e}$, so that
$\sigma (f,e)=1$. Writing
$\rho =ze+\zeta f$, we see that $z,\zeta$ are complex symplectic
coordinates, and $b(\rho )=\mu z\zeta $.

\par Write $e={1\over \sqrt{2}}(e_1+ie_2)$ with $e_1,e_2$ real, so
that $f={1\over \sqrt{2}}(e_2+ie_1)$. Then we see that $e_1,e_2$ is a
real symplectic basis. Also notice that
$$e_1={1\over \sqrt{2}}(e-if),\ e_2={1\over \sqrt{2}}(f-ie),$$
so if $$\rho =ze+\zeta f=xe_1+\xi e_2,$$
we see that $x,\xi $ are real symplectic coordinates on our
symplectic 2-space and
\ekv{{\rm A}.6}
{b(\rho )={\mu \over 2i}(x^2+\xi ^2).}
The resonant action $z\zeta $ becomes ${1\over 2}(x^2+\xi ^2)$ times a
constant factor.

\bigskip

\centerline{\bf \Pa{}. Parameter dependent case.}
\medskip
\par In some applications (for instance when dealing with an energy
dependent monodromy \op{} ([SjZw])) our \fop{} will depend smoothly on
some real parameter $s$, and then it may happen that the non-\res{}
condition is fulfilled for one value of $s$, say for $s=0$ but not
everywhere in any \neigh{} of that point. In this section we show that the
previous results still apply if we require them to hold only to infinite
order \wrt{} $s$ at $s=0$. We do this by checking the earlier
constructions step by step.

\par Let $A=A_s$ be a real symplectic $2n$-matrix depending smoothly on
$s\in{\rm neigh\,}(0,{\bf R})$, such that $A_s$ satisfies the
assumptions of \Prop{} \Sy{}.2, case b. Then $\log A_0$ can be extended to
a smooth family of real matrices $B_s=\log A_s$ with $e^{B_s}=A_s$,
${\strans B}_s+B_s=0$. (The construction of $B=\log A$ can be
reformulated by writing $B=f(A)$, where $f(z)$ is a suitable \hol{} branch
of the logarithm, defined near the spectrum of $A$. We take
$B_s=f(A_s)$ for the same $f$.)

\par Let $\kappa ^s(\rho )$, $s\in{\rm neigh\,}(0,{\bf R})$ be a smooth
family of \canform{}s with
$\kappa ^s(0)=0$ and assume that $\kappa =\kappa ^0$ fulfills the
assumptions of \Th{} \Sy{}.3, so that
\ekv{\Pa{}.1}
{
\kappa ^0(\rho )=\exp H_{p^0}(\rho )+{\cal O}(\rho ^\infty ),
}
where $p^0$ is unique modulo ${\cal O}(\rho ^\infty )$ once its
quadratic part $p_0^0$ has been fixed in accordance with \Prop{}
\Sy{}.2.b. We want to extend $p^0$ to a smooth real-valued family $p^s$,
with
\ekv{\Pa{}.2}
{\kappa ^s(\rho )=\exp H_{p^s}(\rho )+{\cal O}((s,\rho )^\infty ).}
Define $q^s={\cal O}(\rho ^2)$ as in (\Sy{}.23), so that
\ekv{\Pa{}.3}
{(\kappa ^s)^*\partial _s\kappa ^s=H_{q^s},}
and consider the problem analogous to (\Sy{}.24):
\ekv{\Pa{}.4}
{
q^s(\rho )=\int_0^1 \partial _sp^s\circ \exp tH_{p^s}(\rho ) dt+{\cal
O}((s,\rho )^\infty ). }
Putting $s=0$, we get a unique solution $(\partial _sp^s)_{s=0}={\cal
O}(\rho ^2)$ modulo
${\cal O}(\rho ^\infty )$. If we differentiate $k$ times we get
$$\int_0^1(\partial _s^{k+1}p^s)\circ \exp tH_{p^s}(\rho )dt=\partial
_s^kq^s(\rho )+F_k(p^s,..,\partial _s^kp^s,\rho )+{\cal O}((s,\rho
)^\infty ),$$
and if $p^0,..,(\partial _s^kp^s)_{s=0}={\cal O}(\rho ^2)$ already have
been determined, we get $(\partial _s^{k+1}p^s)_{s=0}={\cal O}(\rho ^2)$
from this equation. It is then clear that $(\Pa{}.4)$ has a solution
which is unique mod ${\cal O}((s,\rho )^\infty )$.

\par Let $\widetilde{\kappa }^s=\exp H_{p^s}$. Then
$$(\kappa ^s)^*\partial _s\kappa ^s=(\widetilde{\kappa }^s)^*(\partial
_s\widetilde{\kappa }^s)+{\cal O}((s,\rho )^\infty ),\ \widetilde{\kappa
}^0=\kappa ^0,$$
and as in the proof of \Th{} \Sy{}.3, we see that (\Pa{}.2) holds.

\par We next look at corresponding families of \fop{}s and we start by
extending the equivalence notions of section \Eq{} to the parameter
dependent case. If $V_j\subset {\bf R}^n$ are open \neigh{}s of $0$ and
$I_j\subset {\bf R}$ are open intervals containing $0$, we say that
$v_j\in C^\infty (I_j\times V_j)$, $j=1,2$ are equivalent if they are
equivalent in the sense of section \Eq{} with $V_j$ there replaced by
$I_j\times V_j$. Similarly, we define equivalence for symbols
$a^{(j)}\in S_{{\rm cl}}^0(I_j\times V_j)$ and the corresponding
notion for \pop{}s.

\par Two \canform{}s $\kappa _{j,s}:{\rm neigh\,}(0,{\bf R}^{2n})\to
{\rm neigh\,}(0,{\bf R}^{2n})$ depending \break smoothly on $s\in {\rm
neigh\,}(0,{\bf R})$ with $\kappa _{j,s}(0)=0$, are said to be
equivalent, if $\kappa _{1,s}(\rho )=\kappa _{2,s}(\rho )+{\cal
O}((s,\rho )^\infty )$. (We write $\equiv$ for the parameter version of
equivalence also.) Again $\kappa _{1,s}\equiv \kappa _{2,s}$ iff $\kappa
_{1,s}^{-1}\equiv \kappa _{2,s}^{-1}$.

\par If we parametrize $\kappa _s$ by a \nondeg{} phase $\phi _s
(x,y,\theta )$ depending smoothly on $s$ (and with $(x,y,\theta
)=(0,0,0)$ corresponding to $\kappa _s(0)=0$) then $\kappa _s$ is
equivalent to smooth \fy{} $\widetilde{\kappa }_s$ iff we can
parametrize $\widetilde{\kappa }_s$ by a phase $\widetilde{\phi
}_s(x,y,\theta )$ which is equivalent to $\phi _s(x,y,\theta )$.

\par Consider a \fy{} $U_s=I(a_s,\phi _s)$ of elliptic
\fop{}s as in (\Eq{}.4), associated to a smooth \fy{} of \canform{}s
$\kappa _s$ as above. Assume $\phi _s(0,0,0)=0$. We say that $U_s$ is
equivalent to a second \fy{} $\widetilde{U}_s$, if $\widetilde{U}_s$ has
an associated \fy{} of \canform{}s $\widetilde{\kappa }_s$ with
$\widetilde{\kappa }_s(0)=0$ and we can represent
$\widetilde{U}_s=I(\widetilde{a}_s,\widetilde{\phi }_s)$, with
$\widetilde{\phi }_s(0,0,0)=0$ and with $\widetilde{a}_s\equiv a_s$,
$\widetilde{\phi }_s\equiv \phi _s$ in the sense of families. (We then
have $\widetilde{\kappa }_s\equiv \kappa _s$.)

\par Let $U_s,V_s$ be two families of elliptic \fop{}s as above, with
$U_0=V_0$ and
$U_s^{-1}hD_sU_s\equiv V_s^{-1}hD_sV_s$ (in the sense of families).
Then $U_s\equiv V_s$. In fact, let $\kappa _s,\widetilde{\kappa }_s$ be
the associated \canform{}s and write $U_s^{-1}hD_sU_s=-P_s$, so that
$P_s$ is a smooth \fy{} of \pop{}s of order 0 with real leading symbol
$p_s(\rho )={\cal O}(\rho ^2)$. Then $\kappa _s$ satisfies $\partial
_s\kappa _s(\rho )=(\kappa _s)_*(H_{p_s}(\rho ))$. Similarly
$\partial
_s\widetilde{\kappa }_s(\rho )=(\widetilde{\kappa
}_s)_*(H_{\widetilde{p}_s}(\rho ))$, where $\widetilde{p}_s\equiv
p_s$ and $\widetilde{\kappa }_0=\kappa _0$, so it follows
that $\widetilde{\kappa }_s\equiv \kappa _s$.

\par Without loss of generality, we may assume that $\widetilde{\kappa
}_s=\kappa _s$. Let $W_s$ be some
fixed elliptic family of \fop{}s associated to $\kappa _s$, then
\ekv{\Pa{}.5}
{U_s\equiv A_sW_s,\ V_s\equiv B_s W_s,}
where $A_s$, $B_s$ are smooth families of \pop{}s. We get
$$\eqalign{U_s^{-1}hD_sU_s&\equiv
W_s^{-1}(A_s^{-1}(hD_sA_s))W_s+W_s^{-1}hD_sW_s,\cr V_s^{-1}hD_sV_s&\equiv
W_s^{-1}(B_s^{-1}(hD_sB_s))W_s+W_s^{-1}hD_sW_s.}$$
It follows that
$$A_s^{-1}hD_sA_s\equiv B_s^{-1}hD_sB_s,\ A_0\equiv B_0,$$
and then as in the proof of Lemma \LF{}.1, that $A_s\equiv B_s$ and hence
that
\ekv{\Pa{}.6}
{
U_s\equiv V_s.
}
We can now prove\medskip

\par\noindent \bf \Th{} \Pa{}.1. \it Let $\kappa _s:{\rm neigh\,}(0,{\bf
R}^{2n})\to {\rm neigh\,}(0,{\bf R}^{2n})$, $s\in{\rm neigh\,}(0,{\bf
R})$ be a smooth family of
\canform{}s with $\kappa_s(0)=0$ and assume that $\kappa _0$ satisfies
the assumptions of \Th{} \Sy{}.3. Let $U^s=I(a_s,\phi _s)$ be a
corresponding smooth
\fy{} of elliptic \fop{}s of order 0 with $\phi _s(0,0,0)=0$. Choose
$\mu _j$ and
$p_0=p_0^0$ as there, so that (\Sy{}.32) holds. Then by that theorem and
\Th{}
\LF{}.2, there exists a real-valued smooth function $p^0=p_0^0+{\cal O}(\rho
^3)$ (uniquely determined ${\rm mod\,}({\cal O}(\rho ^\infty ))$), such
that
$\kappa ^0(\rho )=\exp H_{p^0}(\rho )+{\cal O}(\rho ^\infty )$ and a
corresponding \pop{}
$P^0$ with leading symbol
$p^0$, so that $U^0\equiv e^{-iP^0/h}$. ($P^0$ is uniquely determined up
to equivalence and an integer multiple of $2\pi h$.)

$P^0$ can be extended to a smooth family of \pop{}s $P^s$ so that
\ekv{\Pa{}.7}
{
U^s\equiv e^{-iP^s/h}
}
 in the sense of families. Moreover, the family $P^s$ is unique up to
equivalence for families and an integer multipe of $2\pi h$.
The leading symbol $p_s$ satisfies (\Pa{}.2).\rm\medskip

\par\noindent \bf Proof. \rm Let $Q^s\equiv (U^s)^{-1}hD_sU^s$ be the
logarithmic derivative. We first look for $P^s$ solving
\ekv{\Pa{}.8}
{Q^s\equiv -\int_0^1 e^{itP^s/h}(\partial _sP^s)e^{-itP^s/h}dt.}
As in section \LF{}, we first determine $(\partial _sP^s)_{s=0}$. Then we
can write
\ekv{\Pa{}.9}
{
hD_s(e^{-itP^s/h})=R_{t,s}e^{-itP^s/h,}
}
where $R_{t,s}$ is a well-defined smooth \fy{} of 0th order \pop{}s for
$0\le t\le 1$, $s=0$. We can then differentiate (\Pa{}.8) once \wrt{} $s$
and get for $s=0$:
\ekv{\Pa{}.10}
{
\partial _sQ^s\equiv \int_0^1 e^{itP^s/h}[{i\over h}R_{t,s},\partial
_sP^s]e^{-itP^s/h}dt-\int_0^1 e^{itP^s/h}(\partial _s^2P^s)
e^{-itP^s/h}dt. }
>From this we determine $(\partial _s^2P^s)_{s=0}$. Then $(\partial
_sR_{t,s})_{s=0}$ is well-defined in (\Pa{}.9) and we can differentiate
(\Pa{}.10) once more et c and determine $(\partial _s^kP^s)_{s=0}$ for all
$k$. This means that we get a solution of (\Pa{}.8) and we also see
that this solution is unique modulo equivalence for families. From this
we also get the uniqueness of $P^s$ in the \th{}, for if $P^s$ is as in
the \th{}, then it has to satisfy (\Pa{}.8).

\par It remains to show that $P^s$ in (\Pa{}.8) solves (\Pa{}.7). For that, we
put
$$V^s=e^{-itP^s/h},$$
so that by (\Pa{}.8) and the earlier arguments of section \LF{}:
$$Q^s\equiv (V^s)^{-1}(hD_sV^s).$$
Since $Q^s$ is also the log-derivative of the \fy{} $U^s$ and $U^0=V^0$,
we conclude as in (\Pa{}.6), that $U^s\equiv V^s$, and the proof is
complete.\hfill{$\#$}\medskip

\par We end by indicating how to extend the Birkhoff normal form to the
parameter dependent case. Let
\ekv{\Pa{}.11}
{P^{(s)}\sim p^{(s)}(\rho )+hp_1^{(s)}(\rho )+...,\ (s,\rho )\in{\rm
neigh\,}(0,{\bf R}\times {\bf R}^{2n}),}
be smooth in $(s,\rho )$ with $p^{(s)}(\rho )$ real-valued. Asume that
$$p_0^{(0)}(\rho )={1\over 2}\langle (p^{(0)})''(0)\rho ,\rho \rangle $$
satisfies the non-\res{} condition (\BN{}.2). Then in suitable complex
linear symplectic coordinates, we have
\ekv{\Pa{}.12}
{p_0^{(0)}(\rho )=\sum_1^n \mu _j x_j\xi _j,}
and if we allow the coordinates to depend smoothly on $s$, we get
\ekv{\Pa{}.13}
{p_0^{(s)}(\rho )=\sum_1^n \mu _j^{(s)}x_j\xi _j,}
where $\mu _j^{(s)}$ depend smoothly on $s$ and $\mu _j^{(0)}=\mu _j$.
>From the appendix of section \BN{}, it follows that we can find a real
linear \canform{} $\kappa _0^{(s)}$, depending smoothly on $s$ such that
\ekv{\Pa{}.14}
{p_0^{(s)}\circ \kappa _0^{(s)}(\rho )=\sum_1^n \mu _j^{(s)}x_j\xi _j,}
where the coordinates $(x,\xi )$ are now \indep{} of $s$. After composing
$P^{(s)}$ with $\kappa _0^{(s)}$ we can assume that  we have (\Pa{}.13)
with coordintes $x,\xi $ that are \indep{} of $s$. Notice however that
the non-\res{} condition (\BN{}.2) may be violated for $s\ne 0$ \ably{}
close to $0$.

\par We say that a function $r=r_s(x,\xi )\in C^\infty ({\rm
neigh\,}(0,{\bf R}\times {\bf R}^{2n}))$ is resonant if
$H_{p_0^{(0)}}\equiv 0$ in the sense of families. Notice that this
definition does not change, if we replace $H_{p_0^{(0)}}$ by
$H_{p_0^{(s)}}$ (provided we have (\Pa{}.13) in $s$-\indep{}
coordinates). Also notice
$r$ is resonant iff  we have
$r_s\equiv f_s(x_1\xi _1,..,x_n\xi _n)$ for some smooth \fy $f_s$. The
extension of this definition to the case of \pop{}s is immediate. We
next extend Lemma \BN{}.2:
\medskip
\par\noindent \bf Lemma \Pa{}.2. \it For every $v=v^{(s)}\in C^\infty
({\rm neigh\,}(0,{\bf R}^{2n}))$ there exist $u=u^{(s)}$ and
$r=r^{(s)}$ in $C^\infty ({\rm neigh\,}(0,{\bf R}\times {\bf R}^{2n}))$
with $r^{(s)}$ resonant, such that
\ekv{\Pa{}.15}
{H_{p^{(s)}}u^{(s)}=v^{(s)}+r^{(s)}.}
If $v={\cal O}(s^k\rho ^m)$, then we can choose $u,r$ with the same
property.\rm\medskip
\par\noindent \bf Proof. \rm Lemma \BN{}.2 gives a solution
$u^{(0)},r^{(0)}$ for $s=0$. Differentiate (\Pa{}.15) \wrt{} $s$:
\ekv{\Pa{}.16}
{
H_{p^{(s)}}\partial _su^{(s)}=\partial _s v^{(s)}-\{ \partial
_sp^{(s)},u^{(s)}\} +\partial _sr^{(s)}, }
and put $s=0$. Let $(\partial _s u^{(s)})_{s=0}$, $(\partial
_sr^{(s)})_{s=0}$ be the corresponding solutions to this \e{}, given by
Lemma \BN{}.2, then differentiate (\Pa{}.16) et c. In this way, we get
the Taylor series expansion of $u^{(s)}, r^{(s)}$ \wrt{} $s$, and the
lemma follows.\hfill{$\#$}
\medskip
\par\noindent \bf \Prop{} \Pa{}.3. \it There exists a smooth \fy{} of
\canform{}s $\kappa _s:$\hfill{}\break ${\rm neigh\,}(0,{\bf R}^{2n})\to
{\rm neigh\,}(0,{\bf R}^{2n})$ with $\kappa _s(\rho )=\rho +{\cal O}(\rho
^2)$ and
$$p^{(s)}\circ \kappa _s=p_0^{(s)}+r^{(s)},$$
where $r^{(s)}$ is resonant and ${\cal O}(\rho ^3)$.\rm\medskip

\par The proof is essentially identical to that of \Prop{} \BN{}.3. The
treatment of the \op{}s goes through without any changes, and we get
\medskip
\par\noindent \bf \Th{} \Pa{}.4. \it Let $P^{(s)}$ denote also the
$h$-Weyl quantization of the symbol in (\Pa{}.11). Let
$\widetilde{p}_0^{(s)}$ be given by (\Pa{}.13) in the coordinates for
which (\Pa{}.12) holds and let $\widetilde{P}_0^w$ be the corresponding
quantization. Then there exists a smooth \fy{} of elliptic \fop{}s
$V=V^{(s)}$ associated to $\kappa _0^{(s)}\circ \kappa ^{(s)}$ (cf
(\Pa{}.14) and \Prop{} 5.3) such that
$$(V^{(s)})^{-1}P^{(s)}V^{(s)}\equiv P_0^{(s)}+R^{(s)}$$
in the sense of families, where $R^{(s)}$ is a resonant \pop{} of order
$\le 0$ and with leading symbol $={\cal O}(\rho ^3)$. If $P^{(s)}$ is
\sa{}, then we can choose $V^{(s)}$ unitary (microlocally near 0).\rm

\bigskip
\centerline{\bf References.}
\medskip
\par\noindent [Bi] G.D. Birkhoff, \it Dynamical Systems, \rm volume IX.
A.M.S. Colloquium Publications, New York, 1927.
\smallskip
\par\noindent [DiSj] M. Dimassi, J. Sj{\"o}strand, \it Spectral \asy{}s in
the semi-calssical limit, \rm London Math Soc. Lecture Note Ser. 268,
Cambridge Univ. Press, 1999.
\smallskip
\par\noindent [Fr] J.-P. Fran\c coise, \it Propri{\'e}t{\'e}s de
g{\'e}n{\'e}ricit{\'e} des transformations canoniques, \rm pp 216--260 in
Geometric dynamics. Proceedings, Rio de Janeiro, 1981, J. Palis Jr,
editor, Springer LNM 1007.
\smallskip
\par\noindent [Gui] V. Guillemin, \it Wave trace invariants, \rm Duke
Math. J., 83(2)(1996), 287--352.
\smallskip
\par\noindent [Gus] F.G. Gustavsson, \it On constructing formal
integrals of a Hamiltonian system near an equilibrium point, \rm
Astrophys. J., 71(1966), 670--686.
\smallskip
\par\noindent [H{\"o}] L. H{\"o}rmander, \it The analysis of linear partial
differential operators, I--IV, \rm Grundlehren, Springer, 256(1983),
257(1983), 274(1985), 275(1985).\smallskip
\par\noindent [Ia] A. Iantchenko, \it La forme normale de Birkhoff pour
un op{\'e}rateur int{\'e}gral de Fourier, \rm Asymptotic Analysis,
17(1)(1998), 71--92.
\smallskip
\par\noindent [It] H. Ito, \it Integrable symplectic maps and their
Birkhoff normal forms, \rm T{\^o}hoku Math.J., 49(1997), 73--114.
\smallskip
\par\noindent [MeHa] K.R. Meyer, G.R. Hall, \it Introduction to
Hamiltonian dynamical systems and the N-body problem, \rm Applied Math.
Sci. 90, Springer Verlag, 1992.
\smallskip
\par\noindent [Po1] G. Popov, \it Invariant torii, effective stability,
and quasimodes with exponentially small error terms I. Birkhoff normal
forms, \rm Ann. Henri Poincar{\'e} 1(2)(2000), 223--248.
\smallskip
\par\noindent [Po2] G. Popov, \it Invariant torii, effective stability,
and quasimodes with exponentially small error terms II. Quantum Birkhoff
normal forms, \rm Ann. Henri Poincar{\'e} 1(2)(2000), 249--279.
\smallskip
\par\noindent [Sj] J. Sj{\"o}strand, \it Semi-excited states in
non-degenerate potential wells, \rm  Asymptotic Analysis, 6(1992),
29--43.
\smallskip
\par\noindent
[SjZw] J. Sj{\"o}strand, M. Zworski, \it Quantum monodromy and
semi-classical trace formulae, \rm J. Math. Pures et Appl., to appear.
\smallskip
\par\noindent
[St] S. Sternberg, \it Infinite Lie groups and formal aspects of
dynamics, \rm J. of Math. and Mechanics, 10(3)(1961), 451--476.
\smallskip
\par\noindent [Ze1] S. Zelditch, \it Wave invariants at elliptic closed
geodesics, \rm Geom. Funct. Anal., 7(1997), 145--213.
\smallskip
\par\noindent [Ze2] S. Zelditch, \it Wave invariants for non-degenerate
closed geodesics, \rm Geom. Funct. Anal., 8(1998), 179--217.
\smallskip
\par\noindent [Ze3] S. Zelditch, \it Spectral determination of analytic
bi-axisymmetric plane domains, \rm Geom. Funct. Anal., 10(3)(2000),
628--677.
\end